\documentclass[11pt,reqno]{amsart}
\usepackage{amsmath, latexsym, amsfonts, amssymb,amsthm, amscd,epsfig,enumerate}
\usepackage{subfigure}

\advance\hoffset -.75cm

\usepackage{amsfonts}
\usepackage{latexsym}
\usepackage{amsmath}
\usepackage{amssymb}

\newcommand{\R}{\mathbb R}
\newcommand{\Z}{\mathbb Z}
\newcommand{\N}{\mathbb N}
\newcommand{\Prob}{\mathbb P}
\newcommand{\E}{\mathbb E}

\newcommand{\Zd}{{{\mathbb Z}^d}}

\newcommand{\diff}{\rm d}

\newcommand{\eps}{\varepsilon}

\newcommand{\cL}{\mathcal{L}}

\newcommand{\lra}{\longrightarrow}
\newcommand{\ident}{{\mathchoice {\rm 1\mskip-4mu l} {\rm 1\mskip-4mu l}
{\rm 1\mskip-4.5mu l} {\rm 1\mskip-5mu l}}}

\renewcommand{\thesubfigure}{\arabic{subfigure}}
\makeatletter
\renewcommand{\@thesubfigure}{\tiny Figure \thesubfigure: \space}
\renewcommand{\p@subfigure}{}
\makeatother

\newtheorem{teo}{Theorem}[section]
\newtheorem{lem}[teo]{Lemma}
\newtheorem{cor}[teo]{Corollary}
\newtheorem{rem}[teo]{Remark}

\newtheorem{step}{Step}

\begin{document}

\title[Approximating critical parameters of branching random walks]
{Approximating critical parameters \\ of branching random walks}
\author[D.~Bertacchi]{Daniela Bertacchi}
\address{D.~Bertacchi,  Universit\`a di Milano--Bicocca
Dipartimento di Matematica e Applicazioni,
Via Cozzi 53, 20125 Milano, Italy
}
\email{daniela.bertacchi\@@unimib.it}

\author[F.~Zucca]{Fabio Zucca}
\address{F.~Zucca, Dipartimento di Matematica,
Politecnico di Milano,
Piazza Leonardo da Vinci 32, 20133 Milano, Italy.}
\email{fabio.zucca\@@polimi.it}

\date{}

\begin{abstract}
Given a branching random walk on a graph, we consider two kinds of truncations:
either by inhibiting the reproduction outside a subset of vertices or
by allowing at most $m$ particles per vertex.
We investigate the convergence of weak and strong critical parameters
of these truncated branching random walks to the analogous parameters
of the original branching random walk.
As a corollary, we apply our results to the study of the strong
critical parameter of a branching random walk restricted to the
cluster of a Bernoulli bond percolation.
\end{abstract}

%\date{}

\maketitle



\baselineskip .6 cm

\section{Introduction}
\label{sec:intro}
\setcounter{equation}{0}

%*****************************************************************
%\begin{enumerate}
%\item what is a BRW: existence/construction
%\item
%\end{enumerate}
%
%************************************************************

The BRW is a process which serves as a (rough) model for a
population living in a spatially structured environment
(the vertices of a -- possibly oriented -- graph $(X,\mathcal E(X))$),
where each individual lives in a vertex, breeds and dies at
random times and each offspring is placed (randomly) in one of the
neighbouring vertices.
There is no bound on the number of individuals allowed per site.
The vertices may be thought as small ecosystems or squares of
soil (with their proximity connections -- the edges) and
individuals as animals or plants.
Depending on the parameters involved and on the nature of $(X,\mathcal E(X))$,
the population may face almost sure extinction, global survival
(i.e.~with positive probability at any time there will be at least one individual alive)
or local survival (i.e.~with positive probability at arbitrarily large
 times there will be at least one individual alive in a fixed vertex).
These matters have been investigated by several authors (\cite{cf:HuLalley}, \cite{cf:Ligg1}, \cite{cf:Ligg2}, \cite{cf:MadrasSchi},
\cite{cf:PemStac1}, \cite{cf:Stacey03} only to mention a few,
see \cite{cf:Lyons3} for more references).

Let us be more precise as to the definition of the process and of
the environment. The graph $(X,\mathcal E(X))$ is endowed with a
weight function $\mu:X\times X\to [0,+\infty)$ such that
$\mu(x,y)>0$ if and only if $(x,y)\in\mathcal E(X)$ (in which case
we write $x\to y$). We call the couple $(X,\mu)$ a weighted graph.
We require that there exists $K>0$
such that $k(x):=\sum_{y\in X}\mu(x,y)\le
K$ for all $x\in X$ (other conditions will be stated in Section~\ref{sec:definitions}).

Given $\lambda>0$, the branching random walk (BRW$(X)$ or briefly BRW) is the continuous-time
Markov process $\{\eta_t\}_{t\ge0}$, with configuration space $\N^X$,
 where each existing particle at $x$ has an exponential lifespan
of parameter 1 and, during its life, breeds at the arrival times
of a Poisson process of parameter $\lambda k(x)$
and then chooses to send its offspring to $y$
with probability $\mu(x,y)/k(x)$ (note that $(\mu(x,y)/k(x))_{x,y\in X}$ is the transition matrix
of a random walk on $X$). % (for all $y\in X$).
In the literature one usually finds the particular case $k(x)=1$
for all $x\in X$ (i.e.~the breeding rate is constant among locations -- no
place is more fertile than others)
or, sometimes, the case where $\mu={\mathbb I}_{{\mathcal E}(X)}$ (i.e.~the breeding rate
is proportional to the degree and all edges have the same rate).

Two critical parameters are associated to the BRW: the weak (or global) survival
critical parameter $\lambda_w$ and the strong (or local) survival one $\lambda_s$.
They are defined as
\begin{equation}\label{eq:critical}
\begin{split}
\lambda_w&:=\inf\{\lambda>0:\,\Prob^{\delta_{x_0}}\left(\exists t:\eta_t=\underline{0}\right)<1\}\\
\lambda_s&:=\inf\{\lambda>0:\,\Prob^{\delta_{x_0}}\left(\exists \bar t:\eta_t(x_0)=0,\,\forall t\ge\bar t\right)<1\},
  \end{split}
\end{equation}
where $x_0$ is a fixed vertex, $\underline{0}$ is the configuration with no
particles at all sites and $\Prob^{\delta_{x_0}}$ is the law of the process
which starts with one individual in $x_0$. Note that these values
do not depend on the initial configuration, provided that this configuration
is finite (that is, it has only a finite number of individuals), nor on the choice of $x_0$.
See Section \ref{sec:definitions} for a discussion on
the values of $\lambda_w$ and $\lambda_s$.

When $(X,\mu)$ is infinite (and connected),
the BRW is, so to speak, unbounded in two respects:
the environment, since individuals may live at
arbitrarily large distance from their ancestors
(actually $n$-th generation individuals may live
at distance $n$ from the ancestor), and the colonies'
size, since an arbitrarily large number of individuals
may pile up on any vertex.
Hence it is natural to consider ``truncated'' BRWs where
either space or colonies are bounded, and investigate
the relationship between these processes and the BRW.
Indeed, in the literature one often finds problems tackled first
in finite or compact spaces and then reached through
a ``thermodynamical limit'' procedure. One can see easily that  it is possible
to construct the BRW either from the process on finite sets
(spatial truncation) or from the process on infinite space
and a bound on the number of particles per site (particles
truncation). In both cases the truncated process,
for any fixed time $t$, converges almost surely to the BRW.

First we consider ``spatially truncated'' BRWs. We choose
a family of weighted subgraphs $\{(X_n,\mu_n)\}_{n\in\N}$,
%where the sequence of weights $\{\mu_n\}_{n\ge1}$,
such that  $X_n\uparrow
X$, $\mu_n(x,y)\le \mu(x,y)$, and $\mu_n(x,y)
\stackrel{n\to\infty}{\lra}\mu(x,y)$ for all $x,y$.
The process BRW$(X_n)$ can be seen as the BRW$(X)$ with
the constraint that reproductions outside $X_n$ are deleted
and the ones from $x$ to $y$ ($x,y$ in $X_n$) are removed with
probability $1-\mu_n(x,y)/\mu(x,y)$. It is not difficult to see
that for any fixed $t$, as $n$ goes to infinity, the BRW$(X_n)$
converges to the BRW almost surely. Our first result is that
$\lambda_s(X_n)\stackrel{n\to\infty}{\lra}\lambda_s(X)$
%($\lambda_s(X)$
(the latter being the strong survival critical parameter
of the BRW on $(X,\mu)$).
Indeed we prove a slightly more general result
(Theorem~\ref{th:sen2}) which allows us to show that
if $X=\Z^d$ and $X_n$ is the infinite cluster of the
Bernoulli bond percolation of parameter $p_n$,
where $p_n\stackrel{n\to\infty}{\lra}1$ sufficiently fast,
then $\lambda_s(X_n)\stackrel{n\to\infty}{\lra}\lambda_s(X)$
almost surely with respect to the percolation probability space
(Section~\ref{sec:perc}).

Second we  consider BRWs where at most
$m$ individuals per site are allowed (thus taking values in
$\{0,1,\ldots,m\}^X$). We call this process
BRW$_m$ and denote it by $\{\eta_t^m\}_{t\ge0}$.
Note that if $m=1$ we get the contact process (indeed the BRW$_m$
is sometimes referred to as a ``multitype contact process''
-- see for instance~\cite{cf:Neuh}).
It is easily seen that for all fixed $t$ we have
$\eta_t^m\stackrel{m\to\infty}{\lra}\eta_t$ almost surely (see for instance
\cite{cf:PemStac1} where the authors suggest this limit as a way to
contruct the BRW). Clearly, for all $m\ge1$, one may consider the
critical parameters $\lambda_w^m$ and $\lambda_s^m$ defined
as in \eqref{eq:critical} with $\eta_t^m$ in place of $\eta_t$.
One of the main questions we investigate in this paper is whether
$\lambda^m_w\stackrel{m\to\infty}{\lra}\lambda_w$
and $\lambda^m_s\stackrel{m\to\infty}{\lra}\lambda_s$: to our knowledge this
was still unknown even for the case where $X=\Zd$ with $\mu$ transition matrix of the
simple random walk.

Here is a brief outline of the paper.
In Section \ref{sec:definitions} we state the basic terminology and assumptions needed in the sequel.
Section \ref{sec:spatial} is devoted to the spatial approximation of the strong critical parameter
$\lambda_s$ by finite or infinite sets
(see Theorems~\ref{th:sen} and \ref{th:sen2} respectively).
We note that results on the spatial approximation, in the special case when
$X=\Z^d$ and $\mu$ is the transition matrix of the simple random walk, were obtained in
\cite{cf:MountSchin} using a different approach.
In Section~\ref{sec:roadmap} we introduce the technique we use to prove convergence
of the critical parameters of the BRW$_m$.
The technique is essentially a suitable coupling with a supercritical bond
percolation: this kind of comparison has been widely used in
the literature, see for instance \cite{cf:Dur1}, \cite{cf:BGS1}, \cite{cf:Schi1},
\cite{cf:Schi2} and \cite{cf:DurNeu}. Nevertheless the coupling here is quite
tricky, therefore we describe it in four steps which can be adapted to different
graphs.
In Section~\ref{sec:truncateds}
we prove that
$\lambda^m_s$ converges to $\lambda_s$ under some assumptions of self-similarity of the graph (Theorem~\ref{th:main}).
As a corollary, we have $\lambda^m_s\stackrel{m\to\infty}{\lra}\lambda_s$ for $\Zd$ with the simple
random walk. 
The same approach is used
in Section~\ref{sec:truncatedw} to prove the convergence of the sequence
$\lambda^m_w$ to $\lambda_w$ when $X=\Z^d$ (see Theorem~\ref{th:zdrift} and Corollary~\ref{cor:zd},
and Remark~\ref{rem:cartprod} for a slightly more general class of graphs) or
when $X$ is a homogeneous tree (Theorem~\ref{th:tree}).
The results of Section~\ref{sec:spatial} are applied in Section~\ref{sec:perc}
in order to study the strong critical parameter of a BRW restricted to a random
subgraph generated by a Bernoulli bond percolation process.
Section~\ref{sec:open} is devoted to final remarks and open questions.

\section{Terminology and assumptions}
\label{sec:definitions}

In this section we state our assumptions on the
graph $(X,\mu)$; we also recall the description
of the BRW through its generator and the associated semigroup,
and discuss the values of $\lambda_w$ and $\lambda_s$.

Given the (weighted) graph $(X,\mu)$,
 the degree
 of a vertex $x$, $\mathrm{deg}(x)$
is the cardinality of the set $\{y\in X:x\to y%\mu(x,y)>0
\}$; we
require that $(X,\mu)$ is  with bounded geometry, that is
$\sup_{x\in X}\mathrm{deg}(x)<+\infty$.
Moreover we consider $(X,\mu)$ connected, which by our
definition of $\mu$ (recall that $\mu(x,y)>0$ if and only if
$(x,y)\in\mathcal E(X)$) is equivalent to $\mu^{(n)}(x,y)>0$
for some $n=n(x,y)$, where $\mu^{(n)}$ is the $n$-th power of the matrix $\mu$.
When $(\mu(x,y))_{x,y}$ is stochastic (i.e.~$k(x)=1$ for all
$x\in X$), in order to stress this property
 we use the notation $P$, $p(x,y)$ and $p^{(n)}(x,y)$
instead of $\mu$, $\mu(x,y)$ and $\mu^{(n)}(x,y)$.
Define $d(x,y)=\min\{n:\exists \{x_i\}_{i=0}^n, x_0=x,x_n=y,x_i\to x_{i+1}\}$;
note that this is a true metric on $X$ if and only if $(X,\mu)$ is
non oriented.

\begin{figure}
    \centering
    \subfigure[\tiny $X\times Y$ ($X=Y=\Z$).]{\includegraphics[height=4cm]{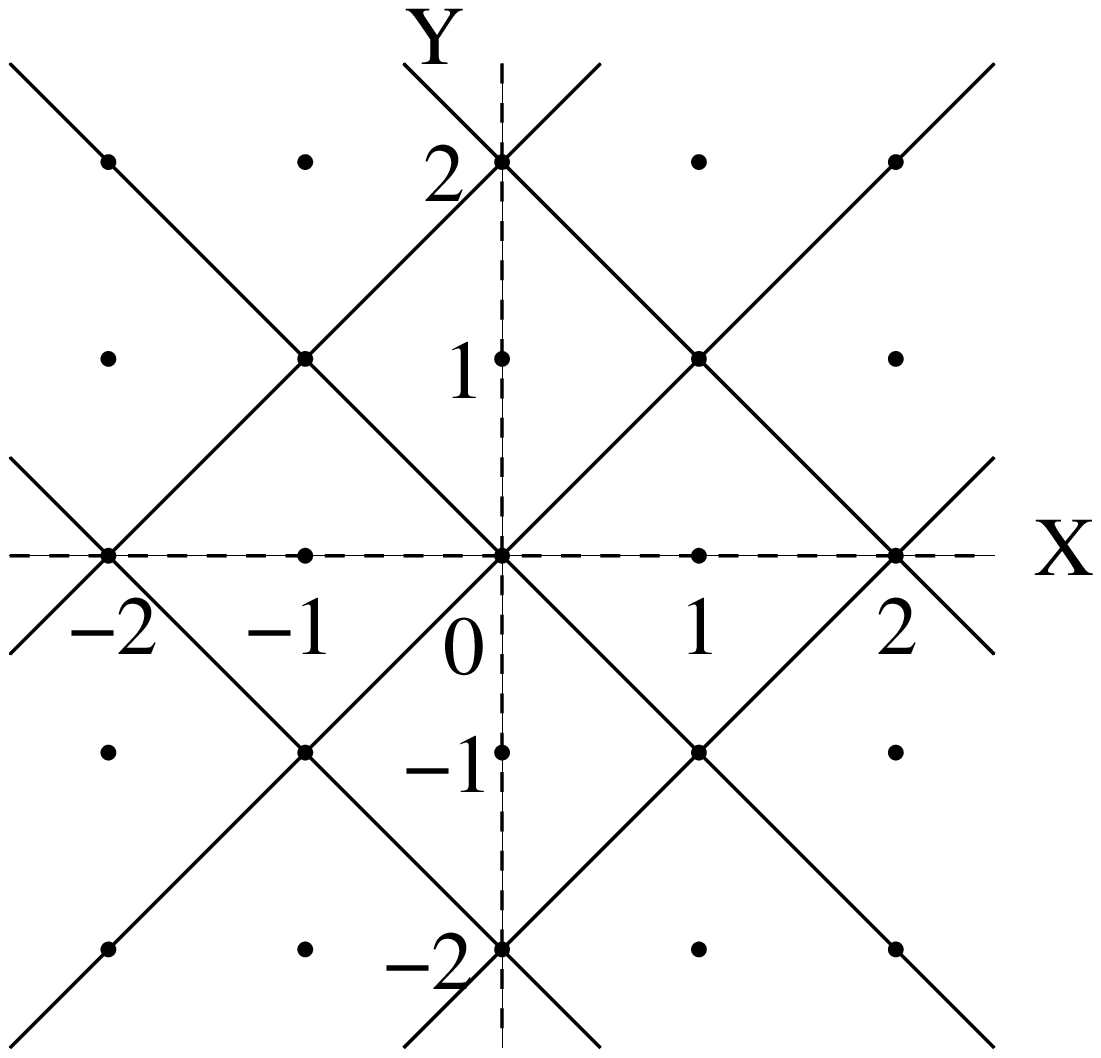}
    \label{fig:xtimesy}}
    \hspace{0.5cm}
    \subfigure[\tiny $X\square Y$ ($X=Y=\Z$).]{\includegraphics[height=4cm]{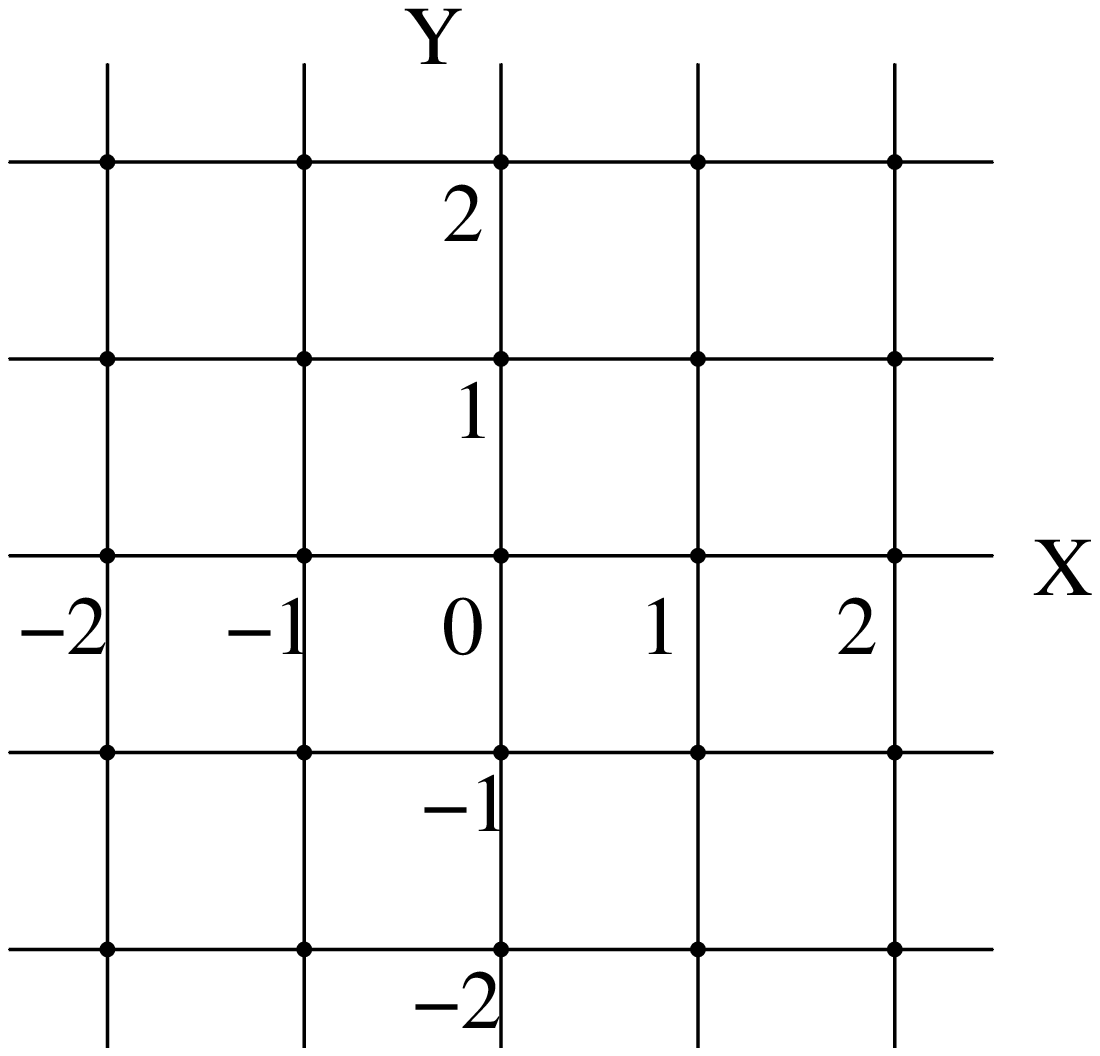}
    \label{fig:xsquarey}}
\end{figure}

We need to define the product of two graphs (in our paper these will
be space/time products):
given two graphs $(X,\mathcal E(X))$, $(Y,\mathcal E(Y))$ we denote by
$(X,\mathcal E(X)) \times (Y,\mathcal E(Y))$ the weighted graph with set of vertices
$X \times Y$ and set of edges $\mathcal E=\{((x,y),(x_1,y_1)): (x,x_1) \in\mathcal E(X),
(y,y_1) \in \mathcal E(Y)\}$ (in  Figure \ref{fig:xtimesy} we draw the connected component of $\Z\times\Z$
containing $(0,0)$).
Besides, by $(X,\mathcal E(X)) \square (Y,\mathcal E(Y))$ we mean the graph
with the same vertex set as before and %set of
vertices
$\mathcal E=\{((x,y),(x_1,y_1)): (x,x_1) \in \mathcal E(X), y=y_1\}
\cup\{ ((x,y),(x_1,y_1)):
x=x_1, (y,y_1) \in \mathcal E(Y)\}$ (see Figure \ref{fig:xsquarey}).

%\begin{figure}
%\centerline{
%\epsfig{figure=xtimesy.eps,height=4cm}}
%\caption{The graph $X\times Y$ ($X=Y=\Z$).}
%\label{fig:xtimesy}
%\end{figure}
%
%\begin{figure}
%\centerline{
%\epsfig{figure=xsquarey.eps,height=4cm}}
%\caption{The graph $X\square Y$ ($X=Y=\Z$).}
%\label{fig:xsquarey}
%\end{figure}

Let $\{\eta_t\}_{t\ge0}$ be the branching random walk on $X$ with parameter
$\lambda$, associated to the weight function $\mu$:
 the configuration space is $\N^X$ and its generator is
%\begin{equation}\label{gengen}
%\cL f (\eta)(z):=  \eta(z)\partial_z^- f(\eta) + \lambda \sum_{x\in X} \eta(x)\,
%\mu(x,z)\,\partial_z^+ f(\eta),
%\end{equation}
\begin{equation}\label{gengen}
\cL f (\eta):=  \sum_{x \in X} \eta(x)\Big (\partial_x^- f(\eta) + \lambda \sum_{y\in X} \,
\mu(x,y)\,\partial_y^+ f(\eta) \Big ),
\end{equation}
where $\partial_x^\pm f(\eta):=f(\eta\pm\delta_x)-f(\eta)$.
%%
%% osservare che non ci sono problemi di esistenza
%%
Analogously the generator of the BRW$_m$
 $\{\eta_t^m\}_{t\ge0}$  is
%\begin{equation}\label{genn}
%\cL_m f (\eta)(z):=  \eta(z)\partial_z^- f + \lambda \sum_{x\in X} \eta(x)\,
%\mu(x,z)\ident_{[0,m-1]}(\eta(z))\,\partial_z^+ f(\eta).
%\end{equation}
\begin{equation}\label{genn}
\cL_m f (\eta):=  \sum_{x \in X} \eta(x)\Big (\partial_x^- f(\eta) + \lambda \sum_{y\in X} \,
\mu(x,y)\ident_{[0,m-1]}(\eta(y))\,\partial_y^+ f(\eta) \Big ),
\end{equation}
Note that the configuration space is still
$\N^X$ (though one may consider $\{0,1,\ldots,m\}^X$ as well).
The semigroup $S_t$ is defined
as
$S_tf(\eta):=\E^\eta(f(\eta_t))$,
where $f$ is any function on $\N^X$ such that the expected value is defined.

The strong and weak survival critical parameters of the BRW clearly
depend on the weighted graph $(X,\mu)$;
we denote them by $\lambda_s(X,\mu)$  and
$\lambda_w(X,\mu)$ (or simply by $\lambda_s(X)$ and $\lambda_w(X)$ or $\lambda_s$ and $\lambda_w$).
 Analogously we denote
by $\lambda^m_s(X,\mu)$ and
$\lambda^m_w(X,\mu)$ (or simply by $\lambda^m_s(X)$ and $\lambda^m_w(X)$ or $\lambda^m_s$ and $\lambda_w^m$) the critical
parameters of the BRW$_m$ on $(X,\mu)$.
It is known (see for instance
\cite{cf:BZ} and \cite{cf:BZ2}) that $\lambda_s=R_\mu:=1/\limsup_n \sqrt[n]{\mu^{(n)}(x,y)}$
(which is easily seen to be independent of $x,y \in X$ since the graph is connected).
On the other hand the explicit value of
$\lambda_w$ is not known in general.
Nevertheless in many cases it is possible to prove that
$\lambda_w=1/\limsup_n  \sqrt[n]{\sum_{y \in X} \mu^{(n)}(x,y)}$
(see \cite{cf:BZ} and \cite{cf:BZ2}).
In particular if $k(x)=K$ for all $x \in X$ then $\lambda_w=1/K$;
thus if $(\mu(x,y))_{x,y}$ is a stochastic matrix then $\lambda_w=1$.

The two critical parameters coincide (i.e.~there is no pure
weak phase) in many cases: if $X$ is finite, or,
when $\mu=P$ is stochastic, if $R=1$.
Here are two sufficient conditions for $R=1$ (when $\mu=P$ is stochastic):
\begin{enumerate}
\item $(X,P)$ is the simple random walk
on a non-oriented graph and the ball of radius $n$ and center $x$ has subexponential growth
 ($\sqrt[n]{|B_n(x)|}\to
1$ as $n \to \infty$).\hfill\break
Indeed for any reversible random walk the following universal lower bound holds
\[
p^{(2n)}(x,x) \ge v(x)/v(B_n(x))
\]
(see \cite[Lemma 6.2]{cf:CGZ}) where $v$ is a reversibility measure. If $P$ is the simple random walk then
$v$ is the counting measure and the claim follows.\hfill\break
An explicit example is the simple random walk on $\Z^d$
or on $d$-dimensional combs (see \cite[Section 2.21]{cf:Woess}
for the definition of comb).
\item
$(X,P)$ is a symmetric, irreducible random walk on an amenable group (see \cite{cf:Woess}).
\end{enumerate}

\section{Spatial approximation}
\label{sec:spatial}
%Throughout this section we suppose that $\sum_{y\in X}\mu(x,y)\le
%k$ then the sum of the convolution powers satisfies $\sum_{y\in
%X}\mu^{(n)}(x,y)\le k^n$.

In this section we consider spatial truncations of the BRW.
We say that $(X,\mu)$ is quasi-transitive
if there exists a %bijection $\gamma$ on $X$ such that the family
%of sets $\{A_x\}_{x\in X}$ where $A_x:=\{\gamma^nx\}_{n\in \N}$
%is finite
%and such that for all $x,y\in X$ we have $p(\gamma(x),\gamma(y))=p(x,y)$,
%.
finite partition of $X$
such that for all couples $(x,y)$ in the same class there exists a bijection $\gamma$ on $X$
satisfying
$\gamma(x)=y$ and, for all  $a,b\in X$, $\mu(\gamma(a),\gamma(b))=\mu(a,b)$
(when the last condition holds we say that $\mu$ is $\gamma$-invariant).
In particular if $\mu(x,y)=p(x,y)$ where $P$ is the simple
random walk on $X$ then it is $\gamma$-invariant
for any automorphism $\gamma$.

In Lemma~\ref{th:sen} and Theorem~\ref{th:sen2} %sequel
$\{X_n\}_{n\in\N}$ will %always
be a sequence of finite
subsets of $X$ such that $X_n\subseteq X_{n+1}$ and
$\bigcup_{n=1}^\infty X_n=X$; we denote by $_n\mu$ the truncation
matrix defined by $_n\mu:=\mu_{|X_n\times X_n}$. We define
$_nR_\mu:=1/\limsup_{k\to\infty}\sqrt[k]{_n\mu^{(k)}(x,y)}$.

\begin{lem}\label{th:sen}
Let $\{X_n\}_{n\in\N}$ be such that $(X_n,\, _n\mu)$ is connected
for all $n$. Then $_nR_\mu\,\ge\,_{n+1}R_\mu$ for all $n$ and
when $X_n\subsetneq X_{n+1}$
we have $_nR_\mu\,>\,_{n+1}R_\mu$. Moreover $_nR_\mu\downarrow R_\mu$.
\end{lem}

\begin{proof}
This is essentially Theorem~6.8 of \cite{cf:Sen}.
\end{proof}

The next result is a generalization of this lemma and it goes beyond the
pure spatial approximation by finite subsets.

\begin{teo}\label{th:sen2}
Let $\{(Y_n,\mu_n)\}_{n\in \N}$ be a sequence of connected weighted
graphs and let $\{X_n\}_{n\in\N}$ be
%a sequence of finite subset of $X$
such that $Y_n \supseteq X_n$.
Let us suppose that
$\mu_n(x,y)\le \mu(x,y)$ for all $n\in\N$, $x,y\in Y_n$
and $\mu_n(x,y)\to\mu(x,y)$ for all $x,y\in X$.
%\hfill\break
If
$(X_n,\, _n\mu)$ is connected for every $n \in \N$
 then $\lambda_s(Y_n,\mu_n)\ge \lambda_s(X,\mu)$ and
$\lambda_s(Y_n,\mu_n)\stackrel{n\to\infty}{
\to} \lambda_s(X,\mu)$.
\end{teo}

\begin{proof}
We note that, for all finite $A\subset X$,
eventually $A\subset Y_m$. Hence
$\mu_n(x,y)$ is well-defined for all sufficiently large $n$.
By Lemma~\ref{th:sen} for any $\eps>0$ there exists $n_0$
such that, for all $n\ge n_0$,
$\lambda_s(X_n,\, _n\mu)=\,_nR_\mu \le R_\mu+\eps/2=\lambda_s(X,\mu)+\eps/2$.
Define $\rho_n={\mu_n}_{|X_{n_0}\times X_{n_0}}$.
Since $X_{n_0}$ is finite and $\rho_n\rightarrow \, _{n_0}\mu$ then
$\lambda_s(X_{n_0},\rho_n)\rightarrow \lambda_s(X_{n_0},\, _{n_0}\mu)$.
Indeed $\lambda_s(X_{n_0},\rho_n)$ and $\lambda_s(X_{n_0},\, _{n_0}\mu)$
are the Perron-Frobenius eigenvalues of $\rho_n$ and $_{n_0}\mu$
respectively and, by construction, for any $\delta>0$, eventually
$(1-\delta)_{n_0}\mu\le\rho_n\le\,_{n_0}\mu$.
If we define $n_1\ge n_0$ such that $\lambda_s(X_{n_0},\rho_{n})<
\lambda_s(X_{n_0},\, _{n_0}\mu)+\eps/2$
for all $n\ge n_1$ then
\[
\lambda_s(Y_{n},\mu_{n})\le \lambda_s(X_{n},{\mu_n}_{|X_{n}\times X_{n}}) \le \lambda_s(X_{n_0},\rho_{n})<
\lambda_s(X_{n_0},\,_{n_0}\mu)+\eps/2<\lambda_s(X,\mu)+\eps,
\]
holds for all $n\ge n_1$.
\end{proof}

A simple situation where the previous theorem applies, is
the non-oriented case ($\mu(x,y)>0$ if and only if $\mu(y,x)>0$)
where $X_n=Y_n$ is the ball of radius $n$ with center at a fixed
vertex $x_0$ of $X$.

\begin{rem}\label{rem:sen}
If $Y_n$ is finite for all $n$, then $\lambda_w(Y_n)=\lambda_s(Y_n)$,
hence $\lambda_w(Y_n)\to\lambda_w(X)$ if and only if $\lambda_w(X)=\lambda_s(X)$.
\end{rem}

%%%%%%%%%%%%%%%%%%%%%%%%%%%%%%%%%%%%%%%%%%%%%%%%%%%%%%%%%%%%%%%%%%
\section{The comparison with an oriented percolation}
\label{sec:roadmap}

From now on, we suppose that $X$ is countable (otherwise $\lambda^n_w=\lambda^n_s=+\infty$).
First of all, we need a coupling between $\{\eta_t\}_{t\ge0}$ and $\{\eta^m_t\}_{t\ge0}$:
think of $\{\eta^m_t\}_{t\ge0}$ as
obtained from $\{\eta_t\}_{t\ge0}$ by removing all the births
which cause more than $m$ particles to live on the same site.
Then we need two other coupled processes.
Fix $n_0\in\N$ and %(the explicit choice of $n_0$ will be stated in Remark \ref{rem:anh}). 
let $\{\bar\eta_t\}_{t\ge0}$ be the process obtained from
the BRW $\{\eta_t\}_{t\ge0}$ by removing all $n$-th generation particles,
with $n> n_0$. Analogously, define $\{\bar\eta^m_t\}_{t\ge0}$
from $\{\eta^m_t\}_{t\ge0}$. Clearly, $\eta_t\ge\bar\eta_t$,
$\eta_t\ge\eta^m_t$, $\eta^m_t\ge\bar\eta^m_t$ and
$\bar\eta_t\ge\bar\eta^m_t$ for all $t\ge0$.
Note that, by construction, the progenies of a given particle
in $\{\bar\eta_t\}_{t\ge0}$ or $\{\bar\eta^m_t\}_{t\ge0}$
lives at a distance from the ancestor not larger than $n_0$
(and the processes go extinct almost surely).

Our proofs of the convergence of $\lambda^m_s$ and $\lambda^m_w$ are essentially divided
in the following four steps.

\begin{step}\label{st:1}
Fix a graph $(I,\mathcal E(I))$ such that the
Bernoulli percolation on $(I,\mathcal E(I))\times\vec\N$
%where $\mathcal V$ is defined as $\{((x,n),(y,n+1)): (x,y)\in\mathcal E(I)\}$
has two phases (where we denote by  $\vec{\N}$ the oriented graph on $\N$, that is, $(i,j)$ is an edge if and only
if $j=i+1$).
\end{step}

\begin{figure}
    \centering
    \subfigure[\tiny $\Z \times \vec{\N}$.]{\includegraphics[height=4cm]{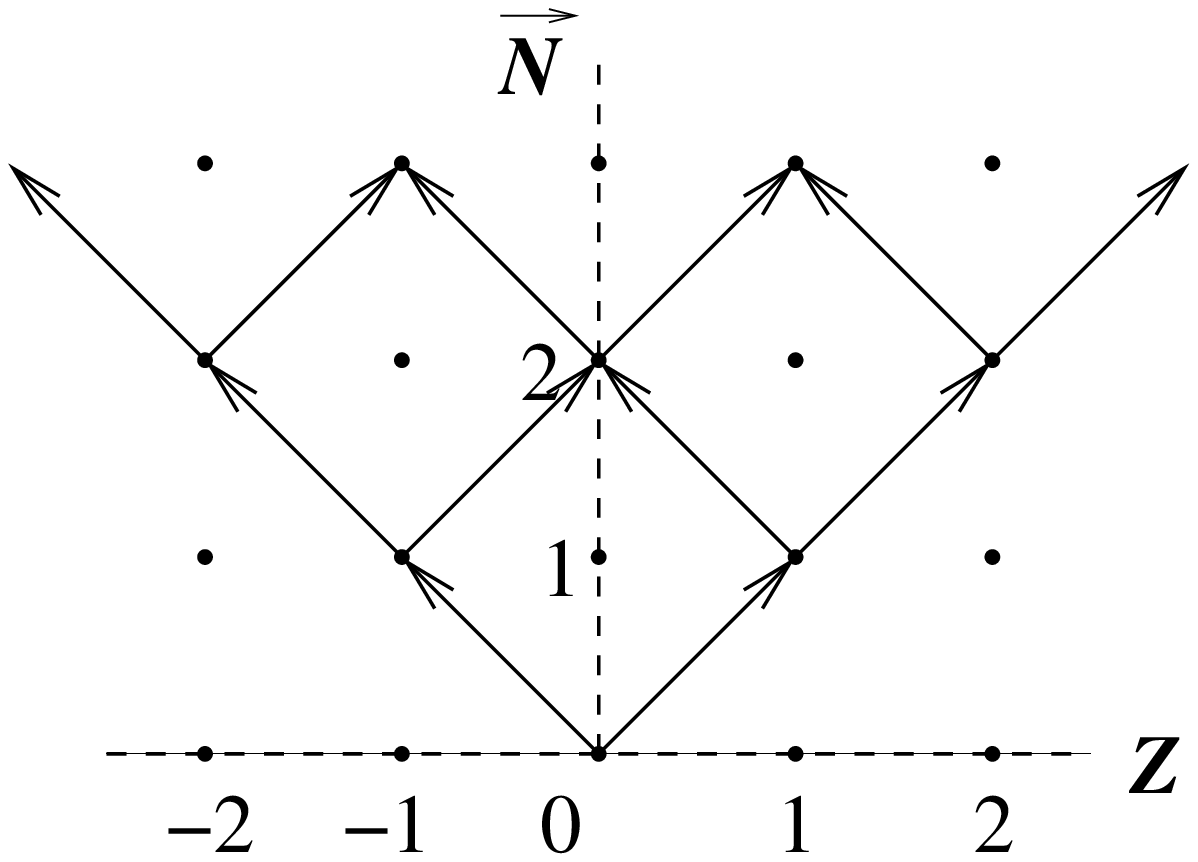}
    \label{fig:ntimesz}}
    \hspace{0.5cm}
    \subfigure[\tiny $\N \times \vec{\N}$.]{\includegraphics[height=4cm]{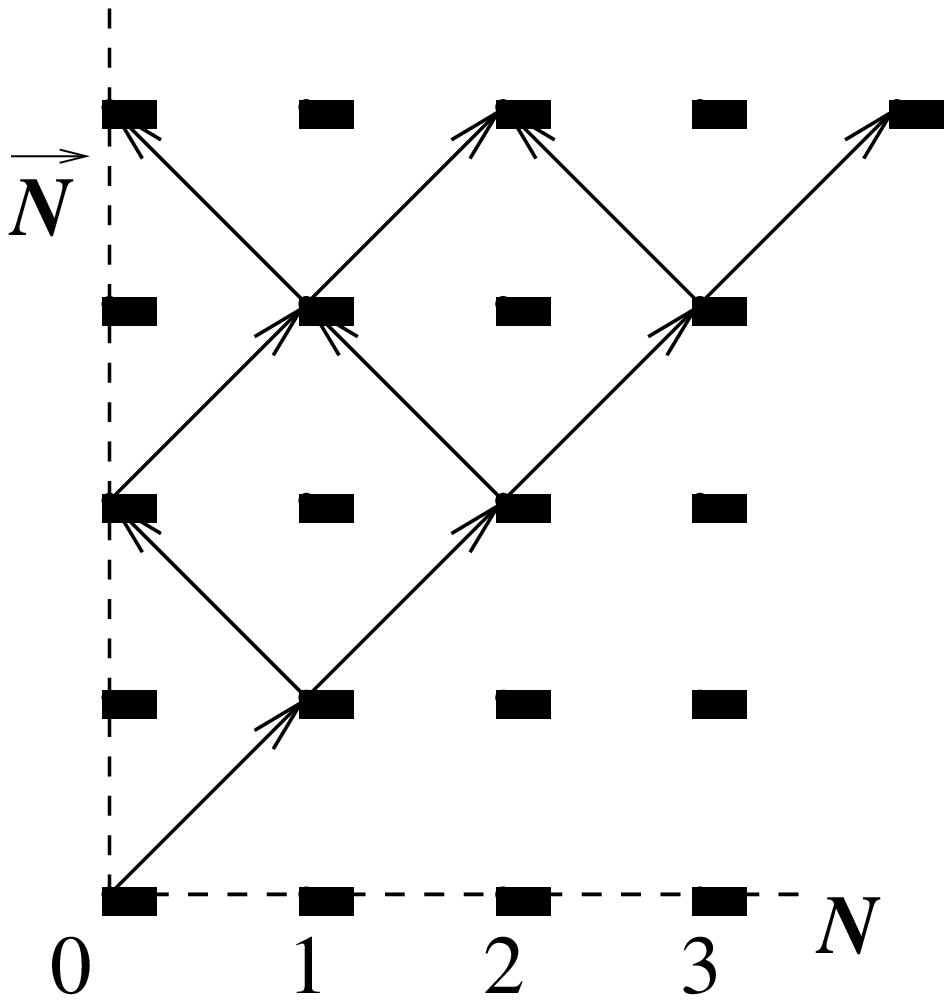}
    \label{fig:ntimesn}}
\end{figure}

%Note that since the Bernoulli bond-percolation on any (non-trivial) cone in $\Z^2$
%has two phases, it is enough to find a copy of the graph
%$\N$ as a subgraph of $I$. This is true for instance
%for any infinite non-oriented graph.
Note that since the (oriented) Bernoulli bond percolation on $\Z \times \vec{\N}$ and
$\N \times \vec{\N}$ has two phases, it is enough to find a copy of the graph
$\Z$ or $\N$ as a subgraph of $I$. This is true for instance
for any infinite non-oriented graph (in this paper, we choose either $I=\Z$, or $I=\N$, or
$I=X$). Figures \ref{fig:ntimesz} and \ref{fig:ntimesn} respectively show the components
of the products $\Z \times \vec{\N}$
and $\N \times \vec{\N}$ containing all the vertices $y$ such that there exists a path from $(0,0)$ to $y$.

\begin{step}\label{st:2}
For all $\lambda>\lambda_w$
(or $\lambda>\lambda_s$ ) and for every $\eps>0$
there exists a collection of disjoint sets $\{A_i\}_{i\in I}$ ($A_i\subset X$ for all $i\in I$), $\bar t>0$,
%$\{t_{i,j}\}_{(i,j)\in\mathcal E(I)}$,
and $k\in\N$, such that, for all $i \in I$,  
%$\forall i \in I$,
%and $t_{i,j}>0$ for any $i\in I$,$\sup_{(i,j)\in\mathcal E(I)}t_{i,j}=\bar t<+\infty$ and for every $i\in I$
\begin{equation}\label{eq:step2}
\Prob \Big ( \forall j:(i,j)\in \mathcal E(I),
\sum_{x\in A_j}\eta_{\bar t%_{i,j}
}(x)\ge k\Big |\eta_{0}= \eta\Big )
>1-\eps,
\end{equation}
for all $\eta$ such that
$\sum_{x\in A_i}\eta(x)= k$ and $\eta(x)=0$ for all $x\not\in A_i$.
The same holds, for some suitable $n_0$, for $\{\bar\eta_t\}_{t\ge0}$ in place of $\{\eta_t\}_{t\ge0}$.
\end{step}

% La richiesta di avere $k$ costante per ogni punto, non e' necessaria. Tuttavia
% appare molto piu` difficile trovare una ragionevole assegnazione a priori
% di una funzione $\{k_x\}_{x \in X}$.

\begin{step}\label{st:3}
Let $\lambda$, $\eps$, $\{A_i\}_{i\in I}$, $\bar t$ and $k$ be chosen as in Step \ref{st:2}.
Then for all sufficiently large $m$, we have that
for all $i \in I$,  
\begin{equation}\label{eq:step3}
\Prob \Big ( \forall j:(i,j)\in \mathcal E(I),
\sum_{x\in A_j}\eta^m_{\bar t%_{i,j}
}(x)\ge k\Big |\eta^m_{0}= \eta\Big )
>1-2\eps,
\end{equation}
for all $\eta$ such that
$\sum_{x\in A_i}\eta(x)= k$,  $\eta(x)=0$ for all $x\not\in A_i$.
The same holds, for some suitable $n_0$, for $\{\bar\eta^m_t\}_{t\ge0}$ in place of $\{\eta^m_t\}_{t\ge0}$.
\end{step}
Step \ref{st:3} is a direct consequence of Step \ref{st:2}.
Indeed let $N_t$ be the total number of particles ever born in the BRW
(starting from the configuration $\eta$) before time $t$; it is clear that $N_t$ is a %continuous time
process bounded above by a %continuos time
branching process with birth rate $K\lambda$, death rate 0 and starting with $k$ particles.
If $N_0<+\infty$ almost surely then for all $t>0$ we have
$N_t<+\infty$ almost surely; hence for all $t>0$ and $\eps>0$
there exists $n(t,\eps)$
such that, for all $i \in I$, 
\[
\Prob\Big (N_t\le n(t,\eps)\Big |\eta_{0}= \eta\Big )>1-\eps,
\]
for all $\eta$ such that
$\sum_{x\in A_i}\eta(x)= k$, $\eta(x)=0$ for all $x\not\in A_i$.
Define $\bar n=n(\bar t,\eps)$. We note that for any event $A$ such that
$\Prob \Big (A \Big |\eta_{0}= \eta \Big)>1-\eps$ we have
\begin{equation}\label{eq:Nt}
\Prob \Big (A \Big |N_t\le \bar n, \eta_{0}= \eta\Big )\ge
\Prob \Big (A, N_t\le \bar n \Big |\eta_{0}= \eta \Big)\ge 1-2\eps.
\end{equation}
Choose %$m\ge2\bar nH$, where $H\in\N$ is the supremum over $x$ of the number of paths of length $n_0$ (chosen in Remark~\ref{rem:anh}) which contain a fixed vertex $x$; $H$ is finite since
%$(X,\mu)$ is with bounded geometry.
$m\ge\bar n$: then $\eta_t=\eta_t^m$ for all $t\le\bar t$ on $\{N_{\bar t}\le\bar n\}
=\bigcap_{t\le\bar t}\{N_t\le\bar n\}$. Thus,
\eqref{eq:step2} and \eqref{eq:Nt} imply \eqref{eq:step3}.
The claim for $\bar\eta^m_t$ is proven analogously.
%As in Fact 2,
%equation~\eqref{eq:step3} still holds for a process
%$\overline{\eta^m_t}$ obtained from $\eta^m_t$ by killing all the
%reproduction trails based on paths of length larger than $n_0$.

\begin{step}\label{st:4}
For all $\lambda>\lambda_w$
(or $\lambda>\lambda_s$ ) and for every $\eps>0$, for all sufficiently large $m$,
there exists a  one-dependent oriented percolation
on $I\times\vec\N$ (with probability $1-2\eps$ of opening all edges)
such that the probability of survival of the BRW$_m$ is larger than the probability
that there exists an infinite cluster containing $({i_0},0)$.
\end{step}
In order to prove Step \ref{st:4} using Step \ref{st:3}, we need another
auxiliary process, namely $\{\widehat\eta_t\}_{t\ge0}$ defined from
$\bar\eta_t$ suppressing all newborns after that the $\bar n$-th particle is born.
If $m\ge\bar n$, then $\widehat\eta_t\le\bar\eta_t^m\le\eta^m_t$ for all $t\ge0$,
and for all $t\le\bar t$, $\widehat\eta_t=\bar\eta_t^m$ on $\{N_{\bar t}\le\bar n\}$.

Consider an edge $((i,n),(j,n+1))$ in $(I,\mathcal E(I))\times\vec\N$:
let it be open if $\eta_t^m$ has at least $k$ individuals in $A_i$ at time $n\bar t$ and
in $A_j$ at time $(n+1)\bar t$.
% it is the same if I open if there is at least $1$ individuals at time $n\bar t$  in $A_i$ and
% in $A_j$ at time $(n+1)\bar t$; in this case I am opening more edges than before
Thus the probability of weak survival of $\eta_t^m$ is bounded from below by
the probability that there exists an infinite cluster containing $({i_0},0)$
in this percolation on $I\times\vec\N$, and, if $A_{i_0}$ is finite,
the probability of strong survival
is bounded from below by the probability that the cluster contains infinitely many points
in $\{({i_0},l):l\in\N\}$ (we suppose that we start with $k$ particles in $A_{i_0}$).
Let $\nu_1$ be the associated percolation measure.
Unfortunately this percolation is neither independent nor one-dependent:
indeed edges $((i,n),(j,n+1))$ and $((i_1,n_1),(j_1,n_1+1))$ can be considered
independent if $n\neq n_1$, but may be dependent if $n=n_1$.
In fact the opening procedure of the edges $((i,n),(j,n+1))$
and $((i_1,n),(j_1,n+1))$ may depend respectively on two different
progenies of particles overlapping on a vertex $x_0$. This may cause
dependence since if in $x_0$ there are already $m$ particles
newborns are not allowed.

To avoid this difficulty we will choose $m$ sufficiently large and
 consider another percolation on $I\times\vec\N$. Let $\widehat\eta_{i_0,t}$
(constructed from $\eta_t$ with the usual removal rules)
start with $k$ particles in $A_{i_0}$: we open all edges $(i_0,0)\to(j,1)$
if $\widehat\eta_{i_0,\bar t}$ has at least $k$ particles in $A_j$.

Let $n=1$: for all $j$ such that $\widehat\eta_{i_0,\bar t}$ has at least $k$ particles in $A_j$
we start again a process $\{\widehat\eta_{j,t}\}_{t\ge0}$ with initial configuration
given by $k$ particles in $A_j$ (chosen among those $\widehat\eta_{i_0,\bar t}$ had there) and zero elsewhere. Note that the $\{\widehat\eta_{j,t}\}_j$
are independent. With a slight abuse of notation we define
$\widehat\eta_t:=\sum_j\widehat\eta_{j,t-\bar t}$ for all $t\in(\bar t,2\bar t]$. 
Choosing $m$ sufficiently large we have that
$\widehat\eta_t\le\bar\eta_t^m\le\eta_t^m$. Indeed it is enough to choose
$m\ge2\bar nH$, where $H\in\N$ is the supremum over $x$ of the number of paths of length $n_0$
 which contain a fixed vertex $x$; $H$ is finite since
$(X,\mu)$ is with bounded geometry.

We iterate the construction for all $n$, obtaining a percolation
$\nu_2$ such that $\nu_1\ge\nu_2$.
Observe that $\nu_2$ is one-dependent since
the set of open edges from $(i,n)$ depends only on
the progenies of the $k$ particles in $A_{i}$ alive in $A_i$ at time $n\bar t$
(hence on $\widehat\eta_{i,t}$, which are independent).
Finally, by equation~\eqref{eq:Nt}, the probability of opening all edges is at least
$1-2\eps$.

\medskip

We note that the trick is to fix a suitable $(I,\mathcal E(I))$ and prove
Step \ref{st:2} for all $\lambda>\lambda_w$: then by Steps \ref{st:4} and \ref{st:1},
for all sufficiently large $m$, the $\lambda$-BRW$_m$ survives with positive probability
and we deduce that $\lambda^m_w\stackrel{m\to\infty}{\lra}\lambda_w$.
On the other hand, to show that $\lambda^m_s\stackrel{m\to\infty}{\lra}\lambda_s$,
we need to prove Step \ref{st:2} with a choice of at least one $A_i$ finite,
say $A_{i_0}$, and $I$ containing a copy of $\Z$ or $\N$
as a subgraph.
Indeed the infinite open cluster in a supercritical Bernoulli bond percolation in
$\Z \times \vec{\N}$ or $\N \times \vec{\N}$ with probability 1
has an infinite intersection with the set $\{(0,n): n \in \N\}$.
As a consequence, in the supercritical case we have, with positive probability,
an infinite open cluster in $\Z \times \vec{\N}$ (resp.~$\N \times \vec{\N}$)
which contains the origin $(0,0)$
and infinite vertices of the set $\{(0,n): n \ge 0\}$.
This (again by Steps \ref{st:3} and \ref{st:4}) implies that, with
positive probability, the $\lambda$-BRW$_m$ starting with $k$ particles in
$A_{i_0}$ has particles alive in $A_{i_0}$ at arbitrarily large times.
Being $A_{i_0}$ finite yields the conclusion.

\begin{rem}\label{rem:roadmap2}
The previous set of steps represents the skeleton of the proofs of Theorems \ref{th:main} and \ref{th:tree}.
In Theorem \ref{th:zdrift} we need a generalization of this approach.
We sketch here the main differences. We choose an oriented graph
$(W, {\mathcal E}(W))$ and a family of subsets of $X$,
$\{A_{(i,n)}\}_{(i,n) \in W}$ such that
\begin{itemize}
\item $W$ is a subset of the set $\Z \times \N$ (note that this is an inclusion between sets not between graphs);
\item for all $n \in \N$ we have that $\{A_{(i,n)}\}_{i:(i,n) \in W}$ is a collection of disjoint subsets of $X$;
\item $(i,n) \to (j,m)$ implies $m=n+1$.
\end{itemize}
The analog of Step \ref{st:2} is the following:
for all sufficiently large $\lambda$
(for instance $\lambda>\lambda_s$ or $\lambda>\lambda_w$) and for every $\eps>0$,
there exists $\bar t>0$ and
$k\in\N$, such that, for all $n \in \N$, $i \in \Z$, and for all $\eta$ such that
$\sum_{x\in A_{(i,n)}}\eta= k$,
\[
\Prob \Big (\forall j:(i,n) \to (j,n+1),
\sum_{x\in A_{(j,n)}}\eta_{(n+1) \bar t%_{i,j}
}(x)\ge k\Big |\eta_{n \bar t}=\eta\Big )
>1-\eps.
\]
Step \ref{st:3} is the same as before and the percolation described in Step \ref{st:4} now concerns
the graph $(W, {\mathcal E}(W))$ (instead of $(I,\mathcal E(I))\times\vec\N$ as it was before).
\end{rem}

\section{Approximation of $\lambda_s$ by $\lambda_s^m$}
\label{sec:truncateds}

We choose the initial configuration
as $\delta_o$ (where $o$ is a fixed vertex in $X$) and we first study the
expected value of the number of individuals in one site at some time,
that is $\E^{\delta_o}(\eta_t(x))$. This is done using the semigroup $S_t$, indeed if we
define the evaluation maps $e_x(\eta):=\eta(x)$ for any $\eta \in \N^X$ and $x \in X$, then
$\E^\eta(\eta_t(x))=S_t e_x(\eta)$.

By standard theorems (see \cite{cf:EK1}, or,
since $\N^X$ is not locally compact, \cite{cf:LiSp} and \cite{cf:BPZ}),
\[
\left. \frac{\diff}{\diff t} S_t e_x \right |_{t=t_0}=S_{t_0} \cL e_x,
\]
from which we deduce
%but from eq \eqref{gengen} \[\cL e_z(\eta)=
% -\eta(z)+ \lambda \sum_{x\in X} \mu(x,z)\eta(x)=-e_{z}(\eta)+
%\lambda \sum_{x\in X} \mu(x,z)e_{x}(\eta) %&z \neq x_0.\\
%\] From the last equation
\begin{equation}\label{eqdiff1}
\frac{\diff}{\diff t}\E^\eta(\eta_t(x))=
% -S_t e_z(\eta)+ \lambda \sum_{x\in X}\mu(x,z)S_t e_{x}(\eta)=
-\E^\eta(\eta_t(x))+
\lambda \sum_{z\in X}\mu(z,x) \E^\eta(\eta_t(z)).
\end{equation}
It is not difficult to verify that
\begin{equation}\label{eq:noimm}
\E^{\delta_{x_0}}(\eta_t(x))=
\sum_{n=0}^\infty \mu^{(n)}(x_0,x)\frac{(\lambda t)^n}{n!}e^{-t}.
%\sum_{n=0}^\infty p^{(n)}(x_0,z)\frac{(\lambda t)^n}{n!}e^{-t}.
\end{equation}

\begin{rem}\label{rem:Ht}
For all $x,x_0 \in X$, for all $\lambda>0$ and $n\in\N$,
\[
\E^{\delta_{x_0}}(\eta_{n}(x))\ge \mu^{(n)}(x_0,x)\frac{\lambda^{n}n^{n}}{n!}e^{-n}
\sim \mu^{(n)}(x_0,x)\frac{\lambda^{n}}{\sqrt{2\pi n}},
\]
and the same inequality holds, if $n_0\ge n$, with $\bar\eta_n$ in place of $\eta_n$.
\end{rem}
Depending on $\lambda$, we may characterize the behaviour of
the expected number of descendants at a fixed site. 
\begin{lem}\label{th:Ht}
%\begin{enumerate}
%\item
%Let $\{a_n\}_{n \in \N}$ be a nonnegative sequence such that $a_{n+m} \ge a_{n}a_{m}$ for all $n,m \in \N$ and
%define $R=1/\lim_{n \to \infty} \sqrt[nr]{a_{n}}$ for some $r \ge 1$. If $\lambda > R$ then, for all nonnegative integer $k$,
%\[
%\lim_{t \to \infty} e^{-t} \sum_{i=0}^\infty a_i \frac{(\lambda t)^{ir+k}}{(i r +k)!}= +\infty.
%\]
%\item
Let us fix $x \in X$. If $\lambda<R_\mu$ then
$\lim_{t\to+\infty}\E^{\delta_{x_0}}(\eta_t(x))=0$;
if $\lambda>R_\mu$ then\break
$\lim_{t\to+\infty}\E^{\delta_{x_0}}(\eta_t(x))=+\infty$.
%\end{enumerate}
\end{lem}

\begin{proof}
Let $\lambda<R_\mu$. For all $\eps>0$ there exists $n_0$ such that $\mu^{(n)}(x_0,x)<1/(R_\mu-\eps)^n$
for all $n\ge n_0$.
If $\eps=(R_\mu-\lambda)/2$ then $\lambda^n\mu^{(n)}(x_0,x)\le \left(\frac{2\lambda}{R_\mu+\lambda}\right)^n$
for all $n\ge n_0$, hence $\E^{\delta_{x_0}}(\eta_t(x))\le Q(t)e^{-t}+e^{-t(R_\mu-\lambda)/(R_\mu+\lambda)}\to0$ as
$t\to\infty$ ($Q$ is a polynomial of degree at most $n_0-1$).

Let $\lambda >R_\mu$.
If $k, r \in \N$ are such that $\mu^{(k)}(x_0,x)>0$ and $\mu^{(ir)}(x,x)>0$ for all $i \in\N$ then
\[
\E^{\delta_{x_0}}(\eta_t(x)) \ge \mu^{(k)}(x_0,x) e^{-t}\sum_{i=0}^\infty \mu^{(ir)}(x,x) (\lambda t )^{ir+k}/(ri+k)!
\]
Let us define $a_n:=\mu^{(nr)}(x,x)$; clearly $a_{n+m} \ge a_n a_m$ and $R_\mu=1/\lim_{n \to \infty} \sqrt[nr]{a_{n}}$
(since $\{a_n\}_{n\in\N}$ is supermultiplicative then the limit exists).

We prove now that for any nonnegative, supermultiplicative sequence
$\{a_n\}_{n \in \N}$,  if $\lambda > R_\mu$ then
\[
\lim_{t \to \infty} e^{-t} \sum_{i=0}^\infty a_i \frac{(\lambda t)^{ir+k}}{(i r +k)!}= +\infty.
\]
Indeed,
let $n_0$ be such that $a_{n_0} \ge 2^{rn_0}(R+\lambda)^{-rn_0}$ and define
$f(t):=e^{-t} \sum_{i=0}^\infty a_i (\lambda t)^{ir+k}/(i r +k)!$.
Clearly
\[
f(t) \ge e^{-t}\left ( \frac{R+\lambda}{2}\right )^k \sum_{i=0}^\infty  \frac{(\lambda^\prime t)^{i n_0 r+k}}{(in_0r+k)!}
\]
where $\lambda^\prime=2\lambda/(R+\lambda)>1$.
For all $i \in \N$ we have that
\[
\frac{(\lambda^\prime t)^{in_0r+k}}{(in_0r+k)!}+ \frac{(\lambda^\prime t)^{in_0r+k+1}}{(in_0r+k+1)!}
+ \cdots + \frac{(\lambda^\prime t)^{in_0r+k+n_0r-1}}{(in_0r+k+n_0r-1)!} \le \frac{(\lambda^\prime t)^{in_0r+k}}{(in_0r+k)!} \frac{(\lambda^\prime t)^{n_0r}-1}{\lambda^\prime t-1}
\]
thus
\[
f(t) \ge e^{-t} \left ( \frac{R+\lambda}{2} \right )^k\frac{\lambda^\prime t-1}{(\lambda^\prime t)^{n_0r}-1} \sum_{i=k}^\infty  \frac{(\lambda^\prime t)^{i}}{i!} \to \infty
\]
exponentially if $t \to \infty$ (since $\lambda^\prime>1$).
\end{proof}
%\begin{lem}\label{th:Ht}
%Let us fix $x \in X$. If $\lambda<R_\mu$ then
%$\lim_{t\to+\infty}\E^{\delta_{x_0}}(\eta_t(x))=0$;
%if $\lambda>R_\mu$ then\break
%$\limsup_{t\to+\infty}\E^{\delta_{x_0}}(\eta_t(x))=+\infty$.
%\end{lem}
%
%\begin{proof}
%Let $\lambda<R_\mu$. For all $\eps>0$ there exists $n_0$ such that $\mu^{(n)}(x_0,x)<1/(R_\mu-\eps)^n$
%for all $n\ge n_0$.
%If $\eps=(R_\mu-\lambda)/2$ then $\lambda^n\mu^{(n)}(x_0,x)\le \left(\frac{2\lambda}{R_\mu+\lambda}\right)^n$
%for all $n\ge n_0$, hence $\E^{\delta_{x_0}}(\eta_t(x))\le Q(t)e^{-t}+e^{-t(R_\mu-\lambda)/(R_\mu+\lambda)}\to0$ as
%$t\to\infty$ ($Q$ is a polynomial of degree at most $n_0-1$).
%
%Let $\lambda >R_\mu$.
%Let $\{n_i\}_{i\in\N}$ be a sequence of natural numbers such that
%$\lim_{i\to\infty}\sqrt[n_i]{\mu^{(n_i)}(x_0,x)}=1/R_\mu$. From this it follows that
%\[
%\E^{\delta_{x_0}}(\eta_{n_i}(x))\ge \mu^{(n_i)}(x_0,x)\frac{\lambda^{n_i}n_i^{n_i}}{n_i!}e^{-n_i}
%\sim \mu^{(n_i)}(x_0,x)\frac{\lambda^{n_i}}{\sqrt{2\pi n_i}},
%\]
%which goes to infinity exponentially as $i\to\infty$ since $\lambda>R_\mu$ and
%eventually $\mu^{(n_i)}(x_0,x)>\left(2/(R+\lambda)\right)^{n_i}$.
%\end{proof}
%%
%% Let $a_n$ be a sequence of nonnegative real numbers such that the corresponding
%% generating function has radius $R$. If $\lambda<R$ then
%% $\lim_{n\to+\infty}\lambda^na_n=0$ and
%% $\lim_{t\to+\infty}e^{-t}\sum_{n} a_n (t\lambda)^t/n!=0$;
%% if $\lambda>R$ then
%% $\limsup_{n\to+\infty}\lambda^na_n=+\infty$ and
%% $\limsup_{t\to+\infty}e^{-t}\sum_{n} a_n (t\lambda)^t/n!=+\infty$.
%%
%%

In the following lemma we prove that, when $\lambda > R_\mu$, if at time 0 we have one
individual at each of $l$ sites $x_1,\ldots,x_l$, then, given any choice
of $l$ sites $y_1,\ldots, y_l$, after some time the expected
number of descendants in $y_i$ of the individual in $x_i$
exceeds 1 for all $i=1,\ldots,l$.

\begin{lem}\label{lem:2.1.5}
Let us consider a finite set of couples
$\{(x_j,y_j)\}_{j=0}^l$; if $\lambda>R_\mu$
then there exists $t=t(\lambda)>0$ %(independent of $M$)
such that
$\E^{\delta_{x_j}}(\eta_t(y_j)) > 1$, % \ge M$,
$\forall j=0,1,\ldots,l$. Moreover, $\E^{\delta_{x_j}}(\bar\eta_t(y_j)) > 1$
when $n_0$ is sufficiently large.
\end{lem}
\begin{proof}
Since $(X,\mu)$ is connected there exist $\{k_j,q_j\}_{j=1,\ldots,l}$
such that, for all $j=1,\ldots,l$ and $n \in \mathbb{N}$,
\[
\mu^{(n+k_j+q_j)}(x_j,y_j) \ge \mu^{(k_j)}(x_j,x_0)\mu^{(n)}(x_0,y_0)\mu^{(q_j)}(y_0,y_j)
\]
and $\mu^{(k_j)}(x_j,x_0)\mu^{(q_j)}(y_0,y_j)>0$.

If $\alpha:=\min_{j=1,\ldots,l}\{\mu^{(k_j)}(x_j,x_0)\mu^{(q_j)}(y_0,y_j)\}$
and $\{n_i\}_{i\in\N}$ is such that %$n_i\ge1$ for all $i$
$\lim_{i \to +\infty} \sqrt[n_i]{\mu^{(n_i)}(x_0,y_0)} =1/R_\mu$,
then for all $j=1,\ldots,l$ and for all $i$
(consider $t=n_i$ and the term with $n=n_i+k_j+q_j$
in the sum \eqref{eq:noimm})
\[
\begin{split}
\E^{\delta_{x_j}}(\eta_{n_i}(y_j)) &\ge \mu^{(n_i)}(x_0,y_0)\frac{(\lambda n_i)^{n_i}}{n_i!}e^{-n_i}\frac{\mu^{(k_j)}(x_j,x_0)\mu^{(q_j)}(y_0,y_j)}{(n_i+k_j+q_j)!}n_i!
(\lambda n_i)^{k_j+q_j}\\
&\ge \mu^{(n_i)}(x_0,y_0)\frac{(\lambda n_i)^{n_i}}{n_i!}e^{-n_i}\frac{\alpha}{(n_i+k_j+q_j)^{k_j+q_j}}
(\lambda n_i)^{k_j+q_j}\\
&\ge \mu^{(n_i)}(x_0,y_0)\frac{(\lambda n_i)^{n_i}}{n_i!}e^{-n_i}\frac{\alpha}{(1+k_j+q_j)^{k_j+q_j}}
\lambda^{k_j+q_j}.
\end{split}
\]
Note that the latter term goes to infinity exponentially as $i \to +\infty$ ($\lambda>R_\mu\ge1$).
Since we have a finite number of sequences, there exists $i_0\ge1$
such that
$\E^{\delta_{x_j}}(\eta_{n_{i_0}}(y_j)) >1$,% \ge M$,
$\forall j=0,1,\ldots,l$.
Choose $t(\lambda)=n_{i_0}$ to conclude. The claim for $\bar\eta$ follows choosing $n_0=n_{i_0}+
\max_j\{k_j+q_j\}$.
\end{proof}

So far we got results on the expected number of individuals, now we show
that, when $\lambda > R_\mu$, for all sufficiently large $k \in \N$, given $k$ particles in a site $x$ at time 0,
``typically'' (i.e.~with arbitrarily large probability)
after some time we will have at least $k$ individuals in each site
of a fixed finite set $Y$.
Analogously, starting with $l$ colonies of size $k$ (in sites
$x_1,\ldots,x_l$ respectively), each of them will, after a
sufficiently long time, spread at least $k$ descendants
in every site of a corresponding (finite) set of sites $Y_i$.

\begin{lem}\label{th:peps} Suppose that $\lambda>R_\mu$.
\begin{enumerate}
\item  %and $\eta_0=k\delta_x$.
Let us fix $x\in X$, $Y$ a finite subset
of $X$ and $\eps>0$.
Then there exists $t=t(\lambda,x)>0$ (independent of $\eps$),
%=t(\eps)$,
$k(\eps,x,Y,\lambda)$ such that, for all $k \ge k(\eps,x,Y,\lambda)$,
\[
\Prob \left(\bigcap_{y \in Y}(\eta_{t}(y)\ge k) \Big | \eta_0(x)=k
\right)>1-\eps.
\]
The claim holds also with $\{\bar\eta_t\}_{t\ge0}$ in place of $\{\eta_t\}_{t\ge0}$ when $n_0$ is
sufficiently large.
\item Let us fix a finite set of vertices $\{x_i\}_{i=1, \ldots,m}$, a collection
of finite sets $\{Y_i\}_{i=1, \ldots,l}$ of vertices of $X$ and $\eps>0$.
Then there exists $t=t(\lambda,\{x_i\},\{Y_i\})$ (independent of $\eps$), $k(\eps,\{x_i\},\{Y_i\},\lambda)$
such that, for all $i=1, \ldots,l$ and $k \ge k(\eps,\{x_i\},\{Y_i\},\lambda)$,
\[
\Prob\left(\bigcap_{y \in Y_i}(\eta_{t}(y)\ge k)
\Big | \eta_{0}(x_i)=k
\right)>1-\eps.
\]
The claim holds also with $\{\bar\eta_t\}_{t\ge0}$ in place of $\{\eta_t\}_{t\ge0}$ when $n_0$ is
sufficiently large.
\end{enumerate}
\end{lem}

\begin{proof} $ $\\
\begin{enumerate}
\item
If we denote by $\{\xi_t\}_t$ the branching process starting from $\xi_0=\delta_x$
then,
by Lemma~\ref{lem:2.1.5},
%such that $x\to y$, by Lemma~\ref{lem:2.1.5},
we can choose $t$ such that $\E^{\delta_x}(\xi_t(y))>1$ for all $y \in Y$.
%Let us note that
We can write $\eta_t(y)=\sum_{j=1}^k\xi_{t,j}(y)$ where
$\xi_{t,j}(y)$ denotes the number of descendants in $y$ of the $j$-th initial particle; note that
$\{\xi_{t,j}(y)\}_{j\in\N}$
is an iid family with $\E(\xi_{t,j}(y))=\E^{\delta_x}(\xi_t(y))$ and $\mathrm{Var}(\xi_{t,j}(y))=:\sigma^2_{t,y}$.
Since $\xi_{t,j}$ is stochastically dominated by a continuous time branching process with birth rate $\lambda$, it is clear that $\sigma^2_{t,y}<+\infty$.
Thus by the Central Limit Theorem, given any $\delta>0$, if $k$ is sufficiently large,
\[
\delta\ge\left|\Prob\left(\sum_{j=1}^k\xi_{t,j}(y)\ge z\right)
-1+\Phi\left(\frac{z-k\E^{\delta_x}(\xi_t(y))}{\sqrt k \sigma_{t,y}}\right)\right|
\]
uniformly with respect to $z\in\R$.
Whence there exists $k(\delta,x,y)$ such that, for all $k \ge k(\delta,x,y)$,
\[
\Prob\left(\eta_{t}(y)\ge k
\right)\ge
1-\Phi\left(\sqrt k\,\frac{1-\E^{\delta_x}(\eta_t(y))}{\sigma_{t,y}}\right)-\delta\ge
1-2\delta,
\]
since $\sqrt k(1-\E^{\delta_x}(\eta_t(y))/\sigma_{t,y}\to -\infty$ as $k\to+\infty$.
Take $k(\delta,x,Y):=\max_{y\in Y}k(\delta,x,y)<+\infty$, and let $D$ be the cardinality of $Y$.
Hence, for all $k \ge k_x$,
\[\Prob \left (\bigcap_{y\in Y}(\eta_{t}(y)\ge k) \Big |\eta_0(x)=k
\right)
\ge 1-2D\delta.
\]
The assertion for $\bar\eta_t$ follows from Lemma \ref{lem:2.1.5}.
\item
Let $\{\xi_t\}_{t \ge 0}$ be as before and choose $t$ such that $\E^{\delta_{x_i}}(\xi_t(y))>1$ for all $y \in Y_i$ and for all $i=1, \ldots l$.
According to (1) above we fix $k_{i}$ such that, for all $k \ge k_{i}$,
\[
\Prob\left(\bigcap_{y\in Y_i}(\eta_{t}(y)\ge k) \Big |\eta_0(x_i)=k
\right)
\ge 1-\eps.
\]
Take $k \ge \max_{i=1, \ldots, l} k_{i}$ to conclude. The assertion for $\bar\eta_t$ follows from Lemma \ref{lem:2.1.5}.
\end{enumerate}
\end{proof}

%What we need to get (2) of the roadmap
% is $\sup_x k_x=k<+\infty$ and
%%$\sup_{x,y} t_{x,y}=\bar t<+\infty$
%%(which holds for $(X,P)$ quasi-transitive).
%$t_{x,y}$ independent of $x$ and $y$.

%\begin{rem}\label{rem:anh}
%Note that Remark~\ref{rem:Ht} and Lemmas~\ref{lem:2.1.5} and \ref{th:peps}
%still hold for a process $\overline{\eta_t}$ obtained
%from $\eta_t$ by killing all the reproduction trails
%based on paths of length larger than a suitable $n_0$
%(for instance choose $n_0=n_{i}$ with $i$ sufficiently large in Remark~\ref{rem:Ht},
%$n_0=n_{i_0}+\max_j\{k_j+q_j\}$ in Lemma~\ref{lem:2.1.5} and
%an analogous choice in Lemma~\ref{th:peps}).
%\end{rem}

\begin{rem}\label{rem:anh}
Note that Remark~\ref{rem:Ht} and Lemmas~\ref{lem:2.1.5} and \ref{th:peps}
can be restated for the process 
$\{\bar{\eta^m_t}\}_{t\ge0}$ if $m$ is sufficiently large.
\end{rem}

\begin{teo} \label{th:main}$ $\\
If at least one of the following conditions holds
\begin{enumerate}
\item  $(X,\mu)$ is %(non-oriented)
% l'ipotesi non-oriented si puo` togliere perche'
% consideriamo p^{(m)} e il grafo e` irriducibile,
% in tal caso si prende come indici X con in piu` gli
% edge "inversi"
quasi-transitive;
\item $(X,\mu)$ is connected and there exists
$\gamma$
bijection on $X$ such that
\begin{enumerate}
\item
$\mu$ is $\gamma$-invariant;
\item for some $x_0\in X$ we have $x_0=\gamma^nx_0$ if and only if $n=0$;
\end{enumerate}
\end{enumerate}
then
\[
\lim_{m\to+\infty}\lambda^m_s=\lambda_s
%=R_\mu
\ge\lim_{m\to+\infty}\lambda^m_w\ge\lambda_w.
\]
Moreover if $\lambda_s=\lambda_w$ then
$\lambda^m_w\downarrow_{m\to+\infty} \lambda_w$.

\end{teo}

\begin{proof}
Remember that $\lambda_s=R_\mu$.
\begin{enumerate}
\item Let us collect one vertex from each orbit into the (finite)
set $\{x_i\}_{i=1, \ldots, l}$ and let $Y_i:=\{y \in X : x_i \to
y\}$. Fix $\lambda >R_\mu$, $I=X$, $\mathcal E(I)=\{(x,y):(x,y)\in \mathcal E(X)\text{ or
}(y,x)\in \mathcal E(X)\}$ and $A_x=\{x\}$. Note that
$(I,\mathcal E(I))$ coincides with $(X,\mathcal E(X))$
if the latter is non oriented.
By these choices, 
Lemma \ref{th:peps} yields Step 2.
%Remark~\ref{rem:percolation} and
To prove that the percolation on $(I,\mathcal E(I))\times \vec{\N}$
has two phases
(that is, $(I,\mathcal E(I))$ is a suitable choice for Step 1) we note that
the existence of the
supercritical phase for the Bernoulli percolation on $X \times
\vec\N$ follows from the fact that the graph $\mathbb{N}$ is a
subgraph of $X$. Moreover in the supercritical Bernoulli percolation on
$\N \times \vec{\N}$ with positive probability the infinite open
cluster contains $(0,0)$ and intersects the $y$-axis infinitely
often. %(see Remark~\ref{rem:percolation}).
Hence by Steps 3 and 4 we have that, for all sufficently large $m$, 
$\lambda_s^m \le \lambda$ and this yields the result.
\item Lemma
\ref{lem:2.1.5} allows us to fix $t$ such that $\E^{\delta_x}(\bar\eta_t(\gamma x))>1$ and 
$\E^{\delta_{\gamma x}}(\bar\eta_t(x))>1$ whence, by
Lemma \ref{th:peps}, for sufficiently large $n_0$,
\[
\Prob\left(\bar\eta_{t}(\gamma x)\ge k \Big | \bar\eta_{0}(x)=k
\right)>1-\eps \quad \text{and} \quad \Prob\left(\bar\eta_{t}(x)\ge k
\Big | \bar\eta_{0}(\gamma x)=k \right)>1-\eps.
\]
This implies
\[
\Prob\left(\bar\eta_{t}(\gamma^{n} x)\ge k \Big |
\bar\eta_{0}(\gamma^{n-1}x)=k \right)>1-\eps \quad \text{and} \quad
\Prob\left(\bar\eta_{t}(\gamma^{n-1} x)\ge k \Big |
\bar\eta_{0}(\gamma^{n} x)=k \right)>1-\eps
\]
for all $n \in \mathbb{Z}$ since $\mu$ is $\gamma$-invariant; one
more time, Steps 3 and 4 %of the roadmap and Remark~\ref{rem:percolation}
yield $\lambda_s^m \le \lambda$ (for sufficiently large $m$) and
the claim (here $I=\mathbb{Z}$ and $A_i=\{\gamma^i x_0\}$).
\end{enumerate}

\end{proof}

%
% Sotto le ipotesi del teorema precedente (almeno quelle del punto (2))
% si ha che se e` vero il seguente claim
% "Se $p$ e` sufficientemente vicino a 1 nella bernoulli percolation di Z^2
% e $r$ e` una retta fissata (es. bisettrice 1o quadrante)
% allora P(#Giant component \cup r=\infty)>0"
% allora $\lambda_s^k \downarrow R=\lambs_s$
%

%%%%%%%%%%%%%%%%%%%%%%%%%%%%%%%%%%%%%%%%%%%%%%%%%%%%%%%%%%%%%%%%%%
\section{Approximation of $\lambda_w$ by $\lambda_w^m$}
\label{sec:truncatedw}

From now on
we set $\mu(x,y)=p(x,y)$ where $P$ is a stochastic matrix. We
stress that in this case $\lambda_w=1$.
We are concerned with the question whether
$\lambda_w^m\downarrow\lambda_w=1$ or not. Under the hypotheses of
Theorem~\ref{th:main}, this is the case when the BRW has no pure
weak phase (i.e.~$R=1$).
The interesting case is $R>1$. Most natural examples are drifting random walks on $\Z^d$
and the simple random walk on homogeneous trees.
In both cases we show that $\lambda^m_w\stackrel{m\to\infty}{\lra}\lambda_w$.
\begin{teo}\label{th:zdrift}
Let $P$ be a random walk on $\Z$ such that $p(i,i+1)=p$, $p(i,i-1)=q$
and $p(i,i)=1-p-q$ for all $i\in\Z$. Then $\lim_{m\to+\infty}\lambda^m_w=1=\lambda_w$.
\end{teo}

\begin{proof}
 We consider $\alpha,\beta\in(0,1)$, $\alpha\le\beta\le(1+\alpha)/2$
and write
\[\begin{split}
p^{(n)}(0,\alpha n)&
=\sum_{i=\alpha n}^{(1+\alpha)n/2}\binom{n}{i,\ \ i-\alpha n,\ \ n-2i+\alpha n}p^{i}
q^{i-\alpha n}(1-p-q)^{n-2i+\alpha n}
\\
&\ge \binom{n}{\beta n,\ \ (\beta-\alpha)n,\ \ (1-2\beta+\alpha)n}p^{\beta n}
q^{(\beta-\alpha)n}(1-p-q)^{(1-2\beta+\alpha)n}\\
&\stackrel{n\to\infty}{\sim} \frac{1}{2\pi n\sqrt{\beta(\beta-\alpha)(1-2\beta+\alpha)}}\,
\left(\frac{p^\beta q^{\beta-\alpha}(1-p-q)^{1-2\beta+\alpha}}{\beta^\beta(\beta-\alpha)^{\beta-\alpha}
(1-2\beta+\alpha)^{1-2\beta+\alpha}}\right)^n.
\end{split}\]
Thus if $\lambda>1$,
$\E^{\delta_0}(\eta_n(\alpha n))$ is bounded from below by a quantity which is asymptotic to
\[
\frac{1}{(2\pi n)^{3/2}\sqrt{\beta(\beta-\alpha)(1-2\beta+\alpha)}}\,
\left(g_\lambda(\alpha,\beta)\right)^n.
\]
where
\[
g_\lambda(\alpha,\beta)=\frac{\lambda p^\beta q^{\beta-\alpha}(1-p-q)^{1-2\beta+\alpha}}{\beta^\beta(\beta-\alpha)^{\beta-\alpha}
(1-2\beta+\alpha)^{1-2\beta+\alpha}}.
\]
Note that $g_\lambda(p-q,p)=\lambda$, thus we may find $\alpha_1<\alpha_2\le\beta_1<\beta_2$
(with $\beta_i\le(1+\alpha_i)/2$, $i=1,2$) such
that $g_\lambda(x,y)>1$, for all $(x,y)\in [\alpha_1,\alpha_2]\times
[\beta_1,\beta_2]$.
By taking $n=\bar n$ sufficiently large one can find three distinct integers
$d_1$, $d_2$ and $d_3$ such that $\alpha_1 n\le d_1<d_2\le\alpha_2n$,
$\beta_1 n\le d_3\le\beta_2n$ and
$g_\lambda(d_l/n,d_3/n)>1$, $l=1,2$.

By reasoning as in Lemma~\ref{th:peps} we have that, for all $\lambda>1$ and $\eps>0$,
there exists $\bar t$, $k=k(\eps,\lambda)$ such that, for all $i\in\Z$, for all $n_0$
sufficiently large,
\[
\Prob \left(\left.\bar\eta_{\bar t}(i+j)\ge k,j=d_1,d_2\right|
\bar\eta_{0}= k\delta_i\right)
>1-\eps.
\]
Since $k$ and $\bar t$ are independent of $i$ we have proven the general version of Step 2
as stated in Remark \ref{rem:roadmap2} (where $W=\{a(d_1,1)+b(d_2,1):a,b\in\N\}$, $A_{(i,n)}=\{i\}$ and $(i,n) \to (j,n+1)$ if and
only if $j-i=d_1$ or $j-i=d_2$).
\end{proof}

In view of Corollary \ref{cor:zd} and Theorem \ref{th:tree} it is
useful to introduce the concept of local isomorphism which allows
to extend some results from $\Z$ to more general graphs. Given
two weighted graphs $(X,\mu)$ and $(I,\nu)$, we say that a map
$f:X \rightarrow I$ is a {\sl local isomorphism} of $X$ on $I$ if for all $x \in
X$ and $i \in I$ we have $\sum_{z \in f^{-1}(i)}
\mu(x,z)=\nu(f(x),i)$.

In this case it is clear that, if we consider the partition of $X$
given by $\{A_i\}_{i\in I}$ where $A_i:=f^{-1}(i)$, we can easily
compute the expected number of particles alive at time $t$ in
$A_i$ starting from a single particle alive in $x$ at time $0$
\begin{equation}\label{eq:locisom}
\sum_{z \in A_i} \E_\mu^{\delta_x}(\eta_t(z))=\E_\nu^{\delta_{f(x)}}(\xi_t(i))
\end{equation}
(where $\{\xi_t\}_{t\ge0}$ is a branching random walk on $(I,\nu)$),
since $\sum_{z \in f^{-1}(i)} \mu^{(n)}(x,z)=\nu^{(n)}(f(x),i)$,
for all $n \in \N$. We note that the latter depends only on
$f(x)$ and $i$. As a consequence $R_\mu \ge R_\nu$.

\begin{cor}\label{cor:zd}
If $P$ is a translation invariant random walk on
$\Z^d$ then $\lim_{k\to+\infty}\lambda^k_w=1=\lambda_w$.
\end{cor}

\begin{proof}
Let $\{Z_n\}_{n\in\N}$ be a realization of the random walk and
$A_i=\{x\in\Z^d:x(1)=i\}$. Note that
\[
\Prob (Z_{n+1}\in A_j|Z_n=w)=\widetilde p(i,j),\quad\forall w\in
A_i,
\]
where $\widetilde P$ is a random walk on $\Z$ with $p=p(0,e_1)$,
$q=p(0,-e_1)$. Using equation \eqref{eq:locisom} and reasoning as in
the proof of the previous theorem, we conclude.
\end{proof}

\begin{rem}\label{rem:cartprod}
The argument of the previous corollary may be applied to a more
general case: let $(Y,Q)$ be a random walk and $(\Z,P)$ be as in
Theorem~\ref{th:zdrift}. Consider $Y\times\Z$ with transition
matrix $\alpha ({\mathbb I}^Y\times P) +(1-\alpha) (Q\times {\mathbb
I}^\Z)$, where $\alpha\in(0,1)$ and by ${\mathbb I}$ we denote the
identity matrix (on the superscripted space). Using the projection
on the second coordinate one proves that
$\lim_{m\to+\infty}\lambda^m_w=1=\lambda_w$.
\end{rem}

%
% ($\Longrightarrow$). Law of large numbers.
% ($\Longleftarrow$). See Remark \ref{rem:local}
%
%

\begin{teo}\label{th:tree}
If $(X,P)$ is the simple random walk on the homogeneous tree of
degree $r$ then $\lim_{m\to+\infty}\lambda^m_w=1=\lambda_w$.
\end{teo}

\begin{proof}
Fix an end $\tau$ in $X$ and a root $o\in X$ and define the map
$h:X\to \Z$ as the usual height (see \cite{cf:Woess} page 129).
Define $A_k=h^{-1}(k)$, $k\in\Z$ (these sets are usually referred
to as horocycles). The projection of the simple random walk on $X$
onto $\Z$ is a random walk with transition matrix $\widetilde  P$ where
$\widetilde  p(a,a+1)=1-1/r$ and $\widetilde  p(a,a-1)=1/r$. Note that for all $x\in X$
\[
\sum_{y\in A_k}p^{(n)}(x,y)=\widetilde p^{(n)}(h(x),k).
\]
%Let us consider a BRW with $\lambda>1$ on $(X,P)$.
%We remark that
%the distribution of the number of descendants in $A_k$
%of an initial particle in $A_h$ alive at time $t$
%depends only on $k-h$ and its mean is
%$\sum_{n=0}^\infty q^{(n)}(0,k-h)(\lambda t)^ne^{-t}/n!$.
%
%Thus it suffices to prove that there exist $d_1,d_2,n\in\N$,
%$d_1\neq d_2$, and $t>0$ such that\break
%$q^{(n)}(0,d_i)(\lambda t)^ne^{-t}/n!>1$ for $i=1,2$.
%
%De Moivre--Stirling's formula (taking $t=n$) yields
%\[
%q^{(n)}(0,\alpha n)\frac{\lambda^n n^n e^{-n}}{n!}
%\sim \frac{1}{\pi n\sqrt{1-\alpha^2}}\left(
%\frac{2\lambda (r-1)^{(1+\alpha)/2}}{r(1+\alpha)^{(1+\alpha)/2}(1-\alpha)^{(1-\alpha)/2}}
%\right)^n
%\]
%Let us define
%\[
%g(\alpha)=\frac{2\lambda (r-1)^{(1+\alpha)/2}}{r(1+\alpha)^{(1+\alpha)/2}(1-\alpha)^{(1-\alpha)/2}}.
%\]
%Note that $g(1-2/r)=\lambda$, thus we may find $\alpha_1<\alpha_2$ such
%that $g(x)>1$, for all $x\in [\alpha_1,\alpha_2]$.
%By taking $n=\bar n$ sufficiently large one can find two distinct integers
%$d_1$ and $d_2$ such that $\alpha_1 n\le d_1<d_2\le\alpha_2n$.
%
%By reasoning as in Lemma~\ref{th:peps} we have that, for all $i\in\N$
%\[
%\Prob \left(\left.\sum_{x\in A_{i+j}}\eta_{\bar n}(x)\ge k,j=d_1,d_2\right|
%\sum_{x\in A_i}\eta_{0}(x)= k\right)
%>1-\eps.
%\]
%Since $k$ and $\bar t:=\bar n$ are independent of $i$ then,
%as in
%(3) of our roadmap, the same holds for some BRW$_m$
%with $m$ sufficiently large.
%We use Bernoulli bond percolation on $\{a(d_1,1)+b(d_2,1):a,b\in\N\}$.

By using equation \eqref{eq:locisom} and reasoning as in
Lemma~\ref{th:peps} and Theorem~\ref{th:zdrift} we have that, for
all $i\in\N$, and some integers $d_1,\, d_2$, 
\[
\Prob \Big ( \sum_{x\in A_{i+j}}\bar\eta_{\bar n}(x)\ge
k,j=d_1,d_2\Big | \bar\eta_{0}(x)= \eta\Big )
>1-\eps,
\]
for all $\eta$
such that $\sum_{x\in A_i}\eta(x)= k$ and $\eta(x)=0$ if $x\not\in A_i$,
and for all $n_0$ sufficiently large.
The claim follows as in Theorem~\ref{th:zdrift}.
\end{proof}

\section{Branching random walks in random environment}
\label{sec:perc}

We use the results of Section~\ref{sec:spatial} to prove some
properties of the BRW in random environment.

Let $(X,\mu)$ be a non-oriented weighted graph.
We consider any subgraph $(Y, {\mathcal E}(Y))$ of $(X,{\mathcal
E}(X))$ as a weighted subgraph with weight function
$\ident_{{\mathcal E}(Y)}\mu$.

Given any $p \in [0,1]$ we consider the Bernoulli bond percolation
on $(X,{\mathcal E}(X))$ and we define the random weighted
subgraphs $(Y^a, {\mathcal E}(Y^a))$ where $Y^a=X$ and ${\mathcal
E}(Y^a)$ is the random set of edges resulting from the percolation
process. We define $\lambda_s(Y^a):=\inf_{A \in \mathcal{A}}
\lambda_s(A)$ where $\mathcal{A}$ is the random collection of all
the connected components of $Y^a$. This corresponds to the
critical (strong) parameter of a BRW where the initial state is
one particle alive at time $0$ in every connected component of $Y^a$.

On the other hand, if there exists a nontrivial critical parameter
$p_c$ for the Bernoulli percolation on $X$ then, for all $p >p_c$
we denote by $(Y^c, {\mathcal E}(Y^c))$ the infinite cluster and
we consider the critical (strong) parameter $\lambda_s(Y^c)$.
Given a sequence $\{p_n\}_{n\in\N}$ such that $p_n \in [0,1]$ for
all $n \in \N$, we consider the sequences $\{Y^a_n\}_{n\in\N}$ and
$\{Y^c_n\}_{n\in\N}$ as the results of independent Bernoulli percolation
processes on $X$ with parameters $\{p_n\}_{n\in\N}$.

Here is the main result; we note that, even when $X=\Z^d$, we do
not require $\mu$ to be the simple random walk.

\begin{teo}
%$\ $ \\
\begin{enumerate}
\item If $\sum_n (1-p_n) < +\infty$ then $\lambda_s(Y^a_n) \to
\lambda_s(X)$ a.s.
\item If $(X,\mu)$ is quasi-transitive  then
$\lambda_s(Y^a)=\lambda_s(X)$ a.s.
\item If $X=\Z^d$
%, $p_n > p_c$ for all $n \in \N$
and $\sum_n (1-p_n) < +\infty$, then
$\lambda_s(Y_n^c) \to \lambda_s(\Z^d)$ a.s.
\item If $X=\Z^d$,
$\mu$ is translation invariant and $p > p_c$ then $\lambda_s(Y^c)
= \lambda_s(\Z^d)$ a.s.

\end{enumerate}
\end{teo}

\begin{proof}
$\ $ \\
\begin{enumerate}
\item By using the Borel-Cantelli Lemma, we have that any finite
connected subgraph of $X$ is eventually contained in a (random)
connected component of $Y^a_n$ almost surely (since
$\sum_n (1-p_n^k) <+\infty$ for all $k \in \N$).
Theorem~\ref{th:sen2} yields the conclusion.
%
% Borel-Cantelli lemma implies that every edge a.s. is eventually open
%http://forums.pcper.com/showthread.php?t=403424
% Lemma. If $\{p_n(i)\}_n$, $i=1,\ldots,k$ is such that $p_n(i) \in [0,1]$ forall $n,i$ and
% $\sum_n(1-p_n(i))<+\infty$, for all $i$
% then $\sum_n (1-\Prod_i p_n(i)) < +\infty$.
% Proof.
% $\alpha_n(i):=1-p_n(i)$ hence $\sum_n \alpha_n(i) < +\infty$ for all $i$.
% \[
% 0 \leq 1-\Prod_i(1-\alpha_n(i)) = \sum_i \alpha_n(i) -\sum_{i_1 <i_2} \alpha_n(i_1) \alpha_n(i_2)+
% \sum_{i_1<i_2<i_3} \alpha_n(i_1) \alpha_n(i_2) \alpha_n(i_3) - \cdots \leq C \sum_i \alpha_n(i)
% \]
% whence
% \[
% \sum_n (1-\Prod_i p_n(i)) \le C \sum_i \alpha_n(i) < +\infty.
% \]
%
\item In this case if we take an infinite orbit $X_0$ then, by
Borel-Cantelli Lemma, for any $m \in \N$, with probability 1,
$Y^a$ contains a ball $B_m$, centered on a vertex $z \in X_0$ and
of radius $m$, with all open edges. Since the critical parameter
of a ball $\lambda_s(B_m)$ does not depend on how we choose its
center in $X_0$, then using Theorem~\ref{th:sen2}, we have that
$\lambda_s(X) \le \lambda_s(Y^a) \le \lambda_s(B_m)\to
\lambda_s(X)$ as $m \to +\infty$.
\item Note that $p_n >p_c$ eventually, hence
$\lambda_s(Y^c_n)$ is well-defined for all sufficiently large $n$.
%We already know that
%$\mu_n(x,y) \to \mu(x,y)$ almost surely, hence
What we need to prove is
that, almost surely, any edge is eventually connected to the
infinite cluster. To this aim we apply the FKG inequality
obtaining that the probability of the event
``the edge $(x,y)$ is open and connected to the infinite cluster $Y_n^c$'' %=A_n
is bigger than $p_n\theta(p_n)$
(where $\theta(p)$ is the probability that a fixed vertex $x$
is contained in the infinite cluster, when each edge is open with probability $p$).
According to Theorem 8.92 of
\cite{cf:Grimmett}, $\theta$ is a differentiable function on
$[0,1]$ hence $1-p\theta(p) \sim (1-p)(1+\theta^\prime(1))$ and
this implies $\sum_n (1-p_n\theta(p_n)) < +\infty$. The
Borel-Cantelli Lemma yields the conclusion.
%
% $q_n :=\pr(A_n)$ then \sum_n(1-q_n) \le \sum_n (1-p_n\theta(p_n))$ and this converges since
% $(1+\theta^\prime(1))\sum_n (1-p_n)$ converges.
%
% More generally:
% (1) $f(a)=b$, $f^\prime_-(a) \not = 0$ implies $b-f(x) \sim f^\prime_-(a)$ when
% $x \to a^-$ hence, if $x_n \to a^-$ then $\sum_n(b-f(x_n))$ converges
% if and only if $\sum_n(a-x_n)$ does.
% (2) $f$ continuous in $(a-\eps,a]$ and $f^\prime(x) \le K$ for all $x \in (a-\eps,a)$
% then |b-f(x)|\le K |a-x| for all $x \in a-\eps,a)$ whence
% if $x_n \in (a-\eps,a)$ eventually, then $\sum_n(b-f(x_n))$ converges
% if and only if $\sum_n(a-x_n)$ does.
%
\item It is tedious but essentially straightforward to prove that,
for any $m \in \N$, with probability 1, $Y^c$ contains an
hypercube $Q_m$ of side-length $m$ with all open edges; as before,
Theorem~\ref{th:sen2} yields the result.
\end{enumerate}
\end{proof}

\section{Final remarks}
\label{sec:open}

At this point the theory of spatial approximation (see Section \ref{sec:spatial}) is
quite complete as far as we are concerned with the basic questions on the convergence of
the critical parameters. Indeed we proved results in this direction (see Theorem
\ref{th:sen2}) for the strong parameter under reasonable assumptions, while the question
on the weak critical parameter, in the pure spatial approximation by finite subsets, %as stated in Theorem
%\ref{th:sen},
is uninteresting (see Remark \ref{rem:sen}). It is possible to further
investigate the convergence of the sequence of weak critical parameters under the hypotheses
of Theorem \ref{th:sen2} by using the characterization
$\lambda_w=1/\limsup_n  \sqrt[n]{\sum_{y \in X} \mu^{(n)}(x,y)}$ which holds in many cases
(see \cite{cf:BZ} and \cite{cf:BZ2} for details).

As for the approximation of the BRW by BRW$_m$s,
we proved that, on quasi-transitive or ``self-similar'' graphs (in the sense of Theorem~\ref{th:main} (2)),
$\lambda_s^m\downarrow\lambda_s$ as ${m\to\infty}$ and,
if there is no weak phase, on $\Z^d$ or on regular trees, $\lambda_w^m\downarrow\lambda_w$
as ${m\to\infty}$.
Here are some natural questions which, as far as we know, are still open:
\begin{itemize}
\item
can one get rid of the hypothesis of quasi-transitivity or
self-similarity in the case concerning the strong critical parameter?
\item
when $\lambda_s>\lambda_w$, is it still true that
$\lambda_w^m\downarrow_{m\to\infty}\lambda_w$,
at least for Cayley graphs or on quasi transitive graphs?
\end{itemize}

\section*{Acknowledgments}
The authors are grateful to Rick Durrett for his invaluable
suggestions.


\begin{thebibliography}{15}

\bibitem{cf:BGS1}
J.~van den Berg, G.~R.~Grimmett, R.~B.~Schinazi,
Dependent random graphs and spatial epidemics,
        Ann.~Appl.~Probab.~\textbf{8} n.~2 (1998), 317-336.

\bibitem{cf:BPZ}
D.~Bertacchi, G.~Posta, F.~Zucca,
        Ecological equilibrium for restrained random walks,
        Ann.~Appl.~Probab.~\textbf{17} n.~4 (2007), 1117-1137.

\bibitem{cf:BZ}
D.~Bertacchi, F.~Zucca,
Critical behaviors and critical values of branching random walks
on multigraphs, J.~Appl.~Probab.~\textbf{45} (2008), 481-497.

\bibitem{cf:BZ2}
D.~Bertacchi, F.~Zucca,
Characterization of the critical values of branching random walks on
weighted graphs through infinite-type branching processes,
ARXIV:0804.0224.

\bibitem{cf:CGZ}
T.~Coulhon, A.~Grigor'yan, F.~Zucca,
The discrete integral maximum principle and its applications,
  Tohoku Math.~J.~\textbf{57}  (2005),  no. 4, 559--587.

\bibitem{cf:Dur1}
R.~Durrett,
Ten lectures on particle systems,
        Springer Lectures Notes in Mathematics  \textbf{1608}, Springer, 1995.

\bibitem{cf:DurNeu}
R.~Durrett, C.~Neuhauser,
Epidemics with recovery in $D=2$,
        Ann.~Appl.~Probab.~\textbf{1} n.~2 (1991), 189-206.

\bibitem{cf:EK1}
S.~N.~Ethier, T.~G.~Kurz, Markov Processes:
characterization and convergence, Wiley Series in Probability and Mathematical
Statistics, John Wiley and Sons INC, New York, 1986.

\bibitem{cf:Grimmett}
        G.~Grimmett, {Percolation},
        Springer-Verlag, Berlin, 1999.

\bibitem{cf:Harris}
T.E.~Harris, A lower bound for the critical probability in a
certain percolation process,
  \textit{Proc.~Cambridge Philos.~Soc.} (1960) \textbf{56},  13--20.

\bibitem{cf:HuLalley}
I.~Hueter, S.P.~Lalley, {Anisotropic branching random walks
on homogeneous trees},  Probab.~Theory Related Fields  \textbf{116},
(2000),  n.1, 57--88.

\bibitem{cf:Ligg1}
T.M.~Liggett,
{Branching random walks and contact processes on homogeneous trees},
Probab.~Theory Related Fields  \textbf{106},  (1996),  n.4, 495--519.

\bibitem{cf:Ligg2}
T.M.~Liggett,
{Branching random walks on finite trees},
Perplexing problems in probability,  315--330, Progr.~Probab.,
\textbf{44}, Birkh\"auser Boston, Boston, MA, 1999.

\bibitem{cf:LiSp}
T.M.~Liggett, F.~Spitzer, Ergodic theorems for coupled random
walks and other systems with locally interacting components,
Z.~Wahrscheinlichkeitstheorie und Verw. Gebiete \textbf{56} n.4
(1981), 443--468.

\bibitem{cf:Lyons3}
R.~Lyons, {Phase transitions on nonamenable graphs.
Probabilistic techniques in equilibrium and nonequilibrium statistical physics},
J.~Math.~Phys. \textbf{41},  (2000),  n.3, 1099--1126.


\bibitem{cf:MadrasSchi}
N.~Madras, R.~Schinazi, {Branching random walks on trees}, Stoch.~Proc.~Appl.
\textbf{42},  (1992),  n.2, 255--267.

\bibitem{cf:MountSchin}
T.~Mountford, R.~Schinazi, A note on branching random walks on
finite sets, J.~Appl.~Probab.~\textbf{42} (2005),  287--294.

\bibitem{cf:Neuh}
C.~Neuhauser, Ergodic theorems for the multitype contact process,
Probab.~Theory Related Fields  \textbf{91}  (1992),  no. 3-4, 467--506.

\bibitem{cf:PemStac1}
R.~Pemantle, A.M.~Stacey, {The branching random walk and
contact process on Galton--Watson and nonhomogeneous trees},
Ann.~Prob.~\textbf{29}, (2001),
 n.4, 1563--1590.

\bibitem{cf:Sen}
        E.~Seneta, {Non-negative matrices and Markov chains},
        Springer Series in Statistics, Springer, New York, 2006.

\bibitem{cf:Schi2}
R.~Schinazi, {On the role of social clusters in the transmission of infectious diseases},
J.~Theoret.~Biol.~\textbf{225},  (2003),  n.1, 59--63.

\bibitem{cf:Schi1}
R.~Schinazi, {Mass extinctions: an alternative to the Allee effects},
Ann.~Appl.~Probab.~\textbf{15},  (2005),  n.1B, 984--991.

\bibitem{cf:Stacey03}
A.M.~Stacey, {Branching random walks on quasi-transitive graphs},
Combin.~Probab.~Comput.~\textbf{12}, (2003), n.3
 345--358.

\bibitem{cf:Woess}
        W.~Woess, {Random walks on infinite graphs and groups},
        Cambridge Tracts in Mathematics, {\textbf 138},
    Cambridge Univ.~Press, 2000.
\end{thebibliography}
\end{document}